\documentclass[british,english]{article}
\usepackage[T1]{fontenc}
\usepackage[latin9]{inputenc}
\usepackage[a4paper]{geometry}
\usepackage{xcolor}
\usepackage{babel}
\usepackage{verbatim}
\usepackage{calc}
\usepackage{dsfont}
\usepackage{amsmath}
\usepackage{amsthm}
\usepackage{amssymb}
\usepackage{graphicx}
\usepackage{wasysym}
\usepackage[unicode=true,pdfusetitle,
 bookmarks=true,bookmarksnumbered=false,bookmarksopen=false,
 breaklinks=false,pdfborder={0 0 1},backref=false,colorlinks=false]
 {hyperref}

\makeatletter
\theoremstyle{plain}
\newtheorem{thm}{\protect\theoremname}
\theoremstyle{definition}
\newtheorem{defn}[thm]{\protect\definitionname}
\theoremstyle{plain}
\newtheorem{prop}[thm]{\protect\propositionname}
\theoremstyle{remark}
\newtheorem*{rem*}{\protect\remarkname}
\theoremstyle{plain}
\newtheorem{conjecture}[thm]{\protect\conjecturename}
\theoremstyle{plain}
\newtheorem{lem}[thm]{\protect\lemmaname}
\theoremstyle{remark}
\newtheorem*{notation*}{\protect\notationname}
\theoremstyle{plain}
\newtheorem{lyxalgorithm}[thm]{\protect\algorithmname}
\theoremstyle{plain}
\newtheorem{cor}[thm]{\protect\corollaryname}

\@ifundefined{date}{}{\date{}}
\usepackage{dsfont}
\usepackage[nottoc,numbib]{tocbibind}
\usepackage{bbm}

\newcommand\blfootnote[1]{%
  \begingroup
  \renewcommand\thefootnote{}\footnote{#1}%
  \addtocounter{footnote}{-1}%
  \endgroup
}

\makeatother

\usepackage[style=alphabetic, maxbibnames=99]{biblatex}
\addto\captionsbritish{\renewcommand{\algorithmname}{Algorithm}}
\addto\captionsbritish{\renewcommand{\conjecturename}{Conjecture}}
\addto\captionsbritish{\renewcommand{\corollaryname}{Corollary}}
\addto\captionsbritish{\renewcommand{\definitionname}{Definition}}
\addto\captionsbritish{\renewcommand{\lemmaname}{Lemma}}
\addto\captionsbritish{\renewcommand{\notationname}{Notation}}
\addto\captionsbritish{\renewcommand{\propositionname}{Proposition}}
\addto\captionsbritish{\renewcommand{\remarkname}{Remark}}
\addto\captionsbritish{\renewcommand{\theoremname}{Theorem}}
\addto\captionsenglish{\renewcommand{\algorithmname}{Algorithm}}
\addto\captionsenglish{\renewcommand{\conjecturename}{Conjecture}}
\addto\captionsenglish{\renewcommand{\corollaryname}{Corollary}}
\addto\captionsenglish{\renewcommand{\definitionname}{Definition}}
\addto\captionsenglish{\renewcommand{\lemmaname}{Lemma}}
\addto\captionsenglish{\renewcommand{\notationname}{Notation}}
\addto\captionsenglish{\renewcommand{\propositionname}{Proposition}}
\addto\captionsenglish{\renewcommand{\remarkname}{Remark}}
\addto\captionsenglish{\renewcommand{\theoremname}{Theorem}}
\providecommand{\algorithmname}{Algorithm}
\providecommand{\conjecturename}{Conjecture}
\providecommand{\corollaryname}{Corollary}
\providecommand{\definitionname}{Definition}
\providecommand{\lemmaname}{Lemma}
\providecommand{\notationname}{Notation}
\providecommand{\propositionname}{Proposition}
\providecommand{\remarkname}{Remark}
\providecommand{\theoremname}{Theorem}

\begin{document}
\global\long\def\X{\mathbb{X}}%
\global\long\def\e{\mathrm{e}}%

\global\long\def\P{\mathbb{P}}%
\global\long\def\Z{\mathbb{Z}}%
\global\long\def\L{\mathbb{L}}%
\global\long\def\B{\mathbb{B}}%
\global\long\def\N{\mathbb{N}}%
\global\long\def\R{\mathbb{R}}%
\global\long\def\I{\mathds{1}}%
\global\long\def\T{\mathds{T}}%
\global\long\def\Nt{N^{(\T)}}%
\global\long\def\Nx{N^{(\X)}}%

\global\long\def\In{\mathtt{In}}%
\global\long\def\Out{\mathtt{Out}}%
\global\long\def\Inv{\mathtt{In}^{[\,\shortparallel\,]}}%
\global\long\def\Outv{\Out^{[\,\shortparallel\,]}}%

\global\long\def\Inh{\mathtt{In}^{[\fallingdotseq]}}%
\global\long\def\Outh{\Out^{[\fallingdotseq]}}%
\global\long\def\alt{\mathrm{alt}}%

\global\long\def\fc{\leadsto_{\mathrm{ffc}}}%
\global\long\def\s{\mathfrak{s}}%
\global\long\def\b{\mathfrak{b}}%
\global\long\def\cD{\mathfrak{d}}%
\global\long\def\cg{\mathfrak{g}}%

\global\long\def\pp{\mathbbm p}%
\global\long\def\kver{\kappa[\shortparallel]}%
\global\long\def\khor{\kappa[\fallingdotseq]}%

\global\long\def\MODEL{\text{LoRaC}}%

\title{A long-range contact process in a random environment}
\author{Benedikt Jahnel \thanks{Technische Universit\"at Braunschweig, Universit\"atsplatz 2, 38106 Braunschweig, Germany, and Weierstrass Institute Berlin, Mohrenstr. 39, 10117 Berlin, Germany,  \texttt{benedikt.jahnel@tu-braunschweig.de}}\qquad{}Anh
Duc Vu \thanks{Weierstrass Institute Berlin, Mohrenstr. 39, 10117 Berlin, Germany,  \texttt{anhduc.vu@wias-berlin.de}}}
\maketitle
\begin{abstract}
We study survival and extinction of a long-range infection process
on a diluted one-dimensional lattice in discrete time. The infection
can spread to distant vertices according to a Pareto distribution,
however spreading is also prohibited at random times. We prove a phase
transition in the recovery parameter via block arguments. This contributes
to a line of research on directed percolation with long-range correlations
in nonstabilizing random environments. 
\end{abstract}

\section*{Introduction}

The contact process is a classical model for the spread of an infection
through a spatially distributed population, where individuals may
spontaneously lose the infection and become susceptible again. First
introduced in \cite{MR0356292}, the model and its multiple generalisations
still attract a tremendous amount of interest coming from a great
variety of fields, see e.g., \cite{MR4452653,MR4572486,MR4594798}
for rather recent contributions and, important in view of this manuscript,
\cite{MR4455874,MR4463081,MR4651886}, where random environments are
considered. Focussing on the discrete-time version on lattices, the
contact process is equivalent to certain models in oriented percolation.
In particular, the key question of survival and extinction of the
infection in the contact process is in one-to-one correspondence to
the existence and absence of an infinite directed path in the associated
percolation model.

The arguably simplest nontrivial undirected percolation model is the
$\Z^{2}$-lattice with either vertices or edges being open with some
probability $p$ independently from each other. The models are then
called site (respectively bond) percolation models and the modeling
idea is usually that of water flowing through open connected components,
i.e.,\ cluster. Now the standard question is whether water can flow
all the way through, i.e.,\ whether the origin lies in an infinite
cluster with positive probability. If so, we are in the socalled supercritical
percolation phase and in the subcritical phase otherwise. In the particular
example just mentioned, the percolation phase transition for the bond
model happens at $p_{c}=1/2$ \cite{MR575895}. 

However, water can only flow in the direction of gravity, so it is
natural to consider directed edges. A simple directed model is the
north-east model on $\Z^{2}$ where connections only form in the north
and east direction introduced in \cite{MR91567}. As pointed out in
\cite{MR757768}, the directed models may have to be handled quite
differently compared to their undirected counterparts. While results
are often similar, the proofs differ greatly.

As mentioned before, we want to consider contact processes, i.e.,\ infections
in space-time rather than the flow of water under gravity. As seen
in the past pandemic, a multitude of different factors influence this
evolution. We want to focus on the following three aspects: range
of infection, sparse environments and lockdowns. More precisely, in
our model we assume that the infection can spread to distant vertices
with polynomial decay in the probability. Additionally, we permanently
remove lattice points via iid Bernoulli random variables, thereby
diluting the lattice. Similarly but now on the time axis, we independently
mark time points at which the transmission of the infection to other
vertices is prohibited. Based on this random environment, we build
our directed bond-site percolation model.

Let us mention that spatial stretches have already been considered
in \cite{MR1112403}. There, a vertex $(t,x)$ is only open with probability
$p(x)\in\{p_{\mathrm{bad}},p_{\mathrm{good}}\}$ where $p(x)$ does
not depend on time. It is shown that survival occurs if $p_{\mathrm{good}}$
occurs sufficiently often and $p_{\mathrm{good}}$ is sufficiently
large. On the other hand, in \cite{MR4442896}, the case of temporal
stretches (on the bonds) has been studied. Here, survival holds even
for any $p_{\mathrm{good}}>p_{c}$ given that $p_{\mathrm{good}}$
occurs sufficiently often, where $p_{c}$ is the critical parameter
for directed bond percolation. The strategy behind both results is
to consider environment groupings and employ a multiscale analysis,
i.e.,\ $\Z^{2}$ is grouped into boxes at different levels and boxes
are combined to form boxes on higher levels. We will follow this general
idea as well and base our construction on \cite{MR2116736} -- which
we have already used in \cite{jahnel2023continuum} and further extend
in this paper -- where percolation of the randomly stretched (undirected)
lattice on $\Z^{2}$ has been proven. Let us note that this result
has recently been refined in \cite{MR4634238} all the way to the
critical parameter $p_{\mathrm{good}}>p_{c}=1/2$.

Simultanously considering temporal and spatial stretches has its own
challenges. For example in \cite{MR4612663}, the authors were able
to link the existence of a nontrivial phase transition on the (undirected)
$\Z^{2}$-lattice to the moments of the stretches. As mentioned there,
their current method only works with one-dimensional stretches. The
problem in our setting is that spreading in space takes time -- time
which might not be available due to lockdowns. We alleviate this issue
by allowing long-range infections. Let us note that considering a
discrete-time process is not a restriction as a simple discretisation
scheme yields also the continuous time case. 

The paper is \foreignlanguage{british}{organised} as follows:
\begin{itemize}
\item In Section \ref{sec:Main}, we introduce the model as well as the
main result, that is, the phase transition of survival and extinction.
We also give the general idea of the proof in Section \ref{subsec:Idea-of-proof}.
\item Section \ref{sec:Proof-skeleton} introduces the core definitions
and lemmas which allow us to prove the main theorem. Details and their
proofs are given in Section \ref{sec:Bands-Labels-Regularity} and
\ref{sec:Proving-Percolation}.
\item Section \ref{sec:Bands-Labels-Regularity} deals solely with the environment
grouping framework while Section \ref{sec:Proving-Percolation} applies
said framework. In particular, this section deals with so called ``drilling''
(Section \ref{subsec:Drilling}) for the multiscale-renormalisation
argument.
\end{itemize}

\section{A long-range contact process ($\protect\MODEL$) \label{sec:Main}}

The model is given as a bond-site percolation model. We consider a
very long street $\Z$ where each $x\in\Z$ represents a location.
Normally, $x$ contains a house with residents (probability $1-q^{(x)}$),
i.e., a potential host for infections. On the other hand, $x$ might
also just be empty (with probability $q^{(x)}$). Now, assume that
there is an infection starting in house $y$. During the day, the
infection might spread to other houses due to people travelling to
other houses. While trips to far-away destinations are rare, they
still happen considerably often via e.g.\ airplanes (probability
$(1+|y-x|)^{-\alpha}$). Each night, all residents of a house recover
with probability $1-p$. In this setting, the survival of an infection
corresponds to a bond-site percolation problem on $\Z\times\Z$ (with
vertices $(t,x)$).

\begin{figure}

\includegraphics[width=1\columnwidth]{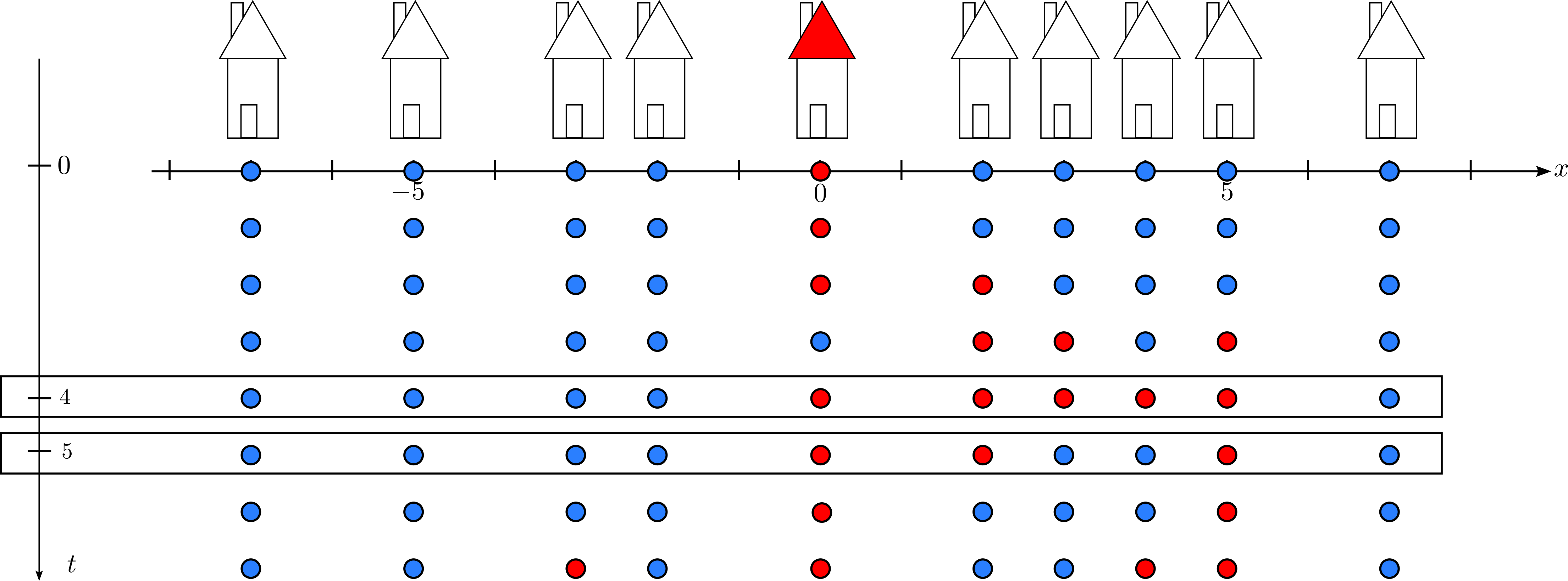}\caption{We start with an infection in the origin which starts infecting other
houses -- preferably close ones. Lockdowns happen at $t=4$ and $5$,
so no spreading occurs during this time. Infections still recover
at any time with probability $1-p$. Note that, here and elsewhere,
we always assume time to flow \textbf{from top to bottom}.\label{fig:houses}}

\end{figure}

During the pandemic, governments have enforced lockdowns during which
people cannot leave their houses. Therefore, no new infections occur
in that time. We mimick this in our model also: Each morning, a global
lockdown is imposed with probability $q^{(t)}$. An illustration of
the model is given in Figure \ref{fig:houses}.

\subsection{Constructing the $\protect\MODEL$ model}

After this verbal discription, let us now give a proper definition
of our model. We highlight that, as mentioned already in the introduction,
contact processes are closely linked to certain directed percolation
problems where the directionality reflects the passing of time.
\begin{defn}[The $\MODEL$]
 Let $q^{(t)},q^{(x)},p\in(0,1)$ as well as $\alpha>1$ be given.
We consider sequences of iid Bernoulli random variables $(\T_{t})_{t\in\Z}$
and $(\X_{x})_{x\in\Z}$ with parameters $q^{(t)}$ respectively $q^{(x)}$.
We call $t$ good if $\T_{t}=1$ and bad otherwise. Analogously, we
call $x$ good if $\X_{x}=1$.

Consider the graph $G=(\Z\times\Z,E)$ where $E$ consists of directed
edges of the form $(t,x)\to(t+1,y)$ with $t,x,y\in\Z$. We study
a mixed bond-site-percolation model on $G$ where all vertices and
edges are open (respectively closed) independently from each other
with probability
\[
\P\{(t,x)\text{ is open}\,\vert\,x\text{ is good}\}=p\qquad\text{and}\qquad\P\{(t,x)\text{ is open}\,\vert\,x\text{ is bad}\}=0\,,
\]
and for an edge $e=\big((t,x)\to(t+1,y)\big)$
\[
\P\{e\text{ is open}\,\vert\,t\text{ is good}\}=(1+|y-x|)^{-\alpha}\qquad\text{and}\qquad\P\{e\text{ is open}\,\vert\,t\text{ is bad}\}=\delta_{xy}
\]
where $\delta_{xy}=1$ iff $x=y$ and $0$ otherwise. We call the
model $\MODEL$ for \emph{long-range contact process}.
\end{defn}

\begin{defn}[Percolation]
 We say that the model \textbf{percolates} if there exists an infinite
sequence of open vertices and edges such that 
\[
(t_{0},x_{0})\to(t_{0}+1,x_{1})\to(t_{0}+2,x_{2})\to\dots
\]
almost surely. In this setting, an infection starting in $x_{0}$
at time $t_{0}$ will spread through open edges and vertices and therefore
survive forever.
\end{defn}

If $\alpha\leq1$, then each vertex has infinitely many outgoing edges
and therefore we already have an infinite number of infected houses
in the first step as well as all subsequent steps. Therefore this
case is trivial. If $\alpha>1$ however, the infection may die out
in certain regimes.
\begin{prop}[Extinction]
 \ 
\begin{enumerate}
\item Let $q^{(t)},q^{(x)}\in(0,1)$ and $\alpha>1$ be given. Then, there
exists $p_{c}\in(0,1)$ such that for every $p<p_{c}$, the model
does not percolate.
\item Let $q^{(t)},q^{(x)},p\in(0,1)$ be given. Then, there exists $\alpha_{c}>1$
such that for every $\alpha>\alpha_{c}$, the model does not percolate.
\end{enumerate}
\end{prop}

\begin{proof}
Point 1 and 2 follow from a simple branching process argument. In
these cases, we completely ignore the environment since it benefits
extinction. $\alpha>1$ implies that the number of potential offsprings
has expectation at most $2\zeta(\alpha)-1$ where
\[
\zeta(\alpha):=\sum_{k=1}^{\infty}k^{-\alpha}\,.
\]
Since each offspring only survives with probability $p$, the actual
number of offsprings is just $(2\zeta(\alpha)-1)\cdot p$, so the
process dies out if we choose $p<(\zeta(\alpha)+1)^{-1}$. (Note that
$\zeta(\alpha)\to1$ as $\alpha\to\infty$.)
\end{proof}
The question then becomes whether survival is actually possible. We
prove a phase transition in the $p$ parameter:
\begin{thm}[Survival via low recovery]
\label{thm:Survival} Let $q^{(t)},q^{(x)}\in(0,1)$ and $\alpha>1$
be given. Then, there exists $p_{c}\in(0,1)$ such that for all $p>p_{c}$,
the $\MODEL$ percolates.
\end{thm}

\begin{rem*}[Continuous time]
 Let us note that this result also holds for the continuous-time
analogue of our model and the proof can be performed via discretisation
arguments.
\end{rem*}
All results also apply for higher dimensions. Survival in $\Z\times\Z$
implies survival in higher dimensions, i.e.,\ $\Z\times\Z^{d}$.
The proof for extinction works analogously as well with $\alpha>d$.

\subsection{Open questions}

Our main theorem is essentially a phase transition in the recovery
of single infections. However, we may also ask ourselves if the process
can survive not by houses staying sick long enough, but rather just
infecting many houses instead. Maybe some clever renormalisation argument
would already do the trick?

\begin{conjecture}[Survival via long spread]
 Given $q^{(t)},q^{(x)},p\in(0,1)$, there exists $\alpha_{c}>1$
such that for every $\alpha\in(1,\alpha_{c})$, the $\MODEL$ percolates.
\end{conjecture}

A different epidemiological concern is the effectiveness of lockdowns
and sparse environments. The comparison of the $\MODEL$ to a Galton--Watson
process with time-dependent offspring distribution tells us that sufficiently
long lockdowns (i.e.\ $q^{(t)}$ close to $1$) will kill off the
infections in the long run. Unfortunately, the effect of the sparse
environment is more complicated to handle.
\begin{conjecture}[Extinction due to sparse environment]
 Given $q^{(t)},p\in(0,1)$ and $\alpha>1$, there exists $q_{c}^{(x)}\in(0,1)$
such that for every $q^{(x)}>q_{c}^{(x)}$, the $\MODEL$ does not
percolate.
\end{conjecture}

We see that infinitely long edges are definitely required for the
model to percolate. If the length of the edges was bounded, then the
whole infection would be confined to a finite region since the infection
is not able to cross over large gaps. However, the exact asymptotic
decay of the edges is crucial and we are currently unable to deal
with the case of exponential decay.
\begin{conjecture}[Fewer edges]
 The $\MODEL$ has a phase transition even if edges are only present
with probability
\[
\exp(-\alpha|y-x|)\,.
\]
\end{conjecture}

This case would be related to the actual ``randomly stretched directed
lattice'' with stretches in both the temporal \cite{MR4442896} and
spatial component \cite{MR1112403}. 

Unfortunately, both ideas cannot be directly combined to prove percolation.
In \cite{MR1112403}, one considers extremely thin boxes where the
height is an exponential of the width. While the multiscale estimates
would still work, the frameworks in \cite{MR2116736,MR4442896} restrict
ourselves to boxes which do not permit the same extreme scaling.

\subsection{Idea of proof \label{subsec:Idea-of-proof}}

\begin{figure}

\includegraphics[width=1\columnwidth]{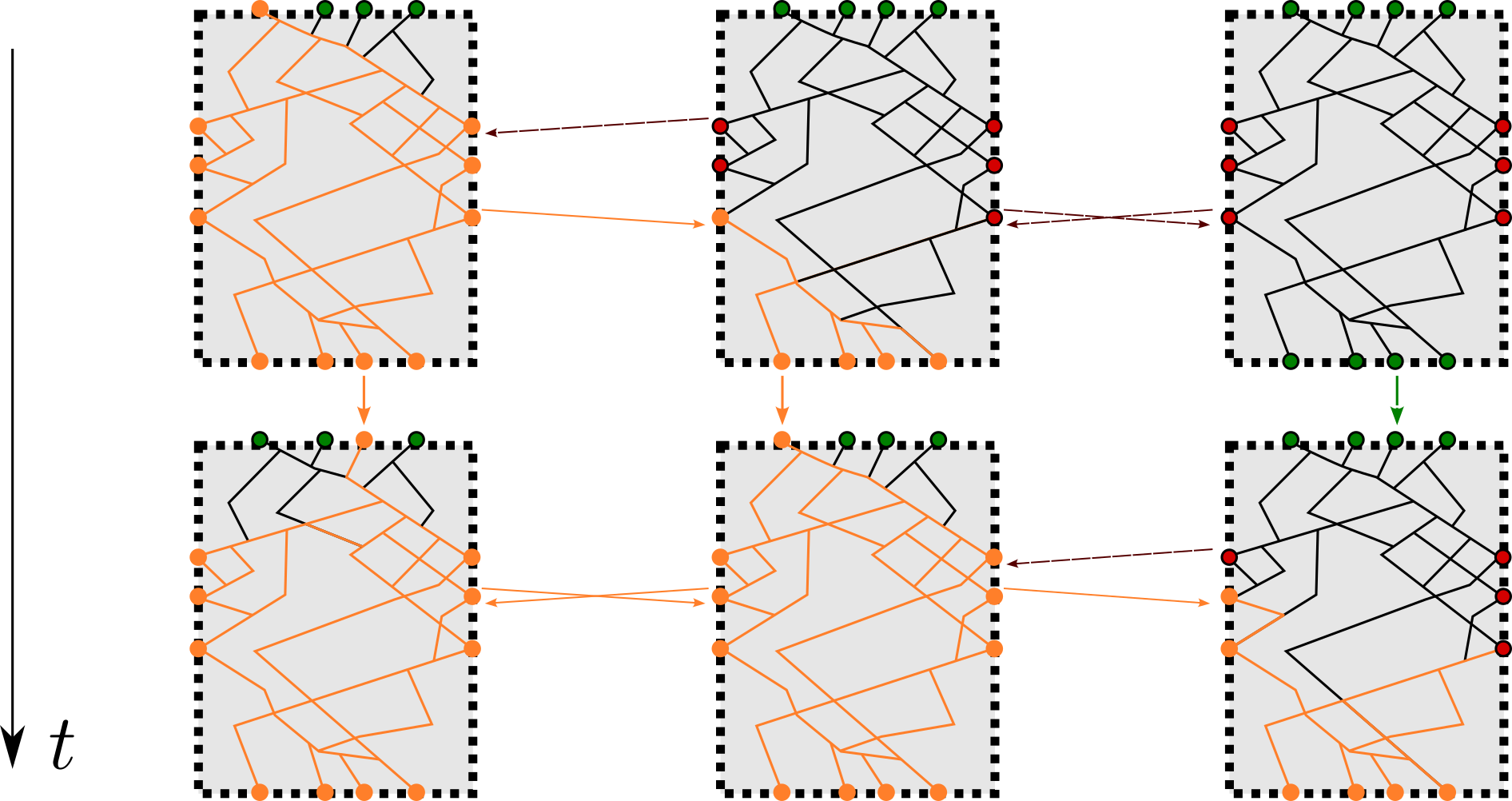}\caption{The environment is divided into boxes at different levels. An orange
vertex starts infecting everything on its way down. Boxes are well
connected since $p$ is large. The infection uses special vertices
(outputs/inputs depicted as circles) to spread to other neighbouring
boxes. The environment between boxes is hostile, so usually only few
connections are found. \label{fig:strategy} }

\end{figure}

The setup for the proof of Theorem \ref{thm:Survival} is quite long
and it is easy to get lost in details. While -- as always -- the
main difficulty lies in those details, they are not as insightful
to the general idea and have already been dealt with in other works.
We will not reinvent the wheel, but building a cart from it has merit
in itself. The procedure is as follows:
\begin{enumerate}
\item We move away from Bernoulli random variables in the $\MODEL$ and
use geometric ones instead. Both model formulations are equivalent
in terms of percolation, but the latter is much more convenient to
use.
\item The next step lies in dividing both the time and space random environments
into bands.
\item From there, we will use these bands to define $n$ boxes: rectangular
subsets in $\Z\times\Z$. These boxes are roughly exponentially large
in $n$ and consist of $n-1$ boxes.
\item Each $n$ box has some special vertices on the boundary which we will
call (horizontal/vertical) inputs and outputs. There are exponentially
many of those vertices.
\item With high probability, $n$ boxes are ``good'' which means that
the aforementioned inputs and outputs are well connected. Also with
high probability, the output of an $n$ box will connect to the input
of a neighbouring $n$ box (restricted by directionality). This is
graphically represented in Figure \ref{fig:strategy}.
\item As $n\to\infty$, the $n$ boxes will always be good which yields
an infinite cluster.
\end{enumerate}
We make this procedure rigorous in the next section.

\section{Proof skeleton \label{sec:Proof-skeleton}}

In the following, we will give the bare proof skeleton leading up
to the main result of phase transition. We try to keep the main ideas
while omitting most details and proofs.

\subsection{Alternative model construction and coupling}

\begin{figure}

\includegraphics[width=1\columnwidth]{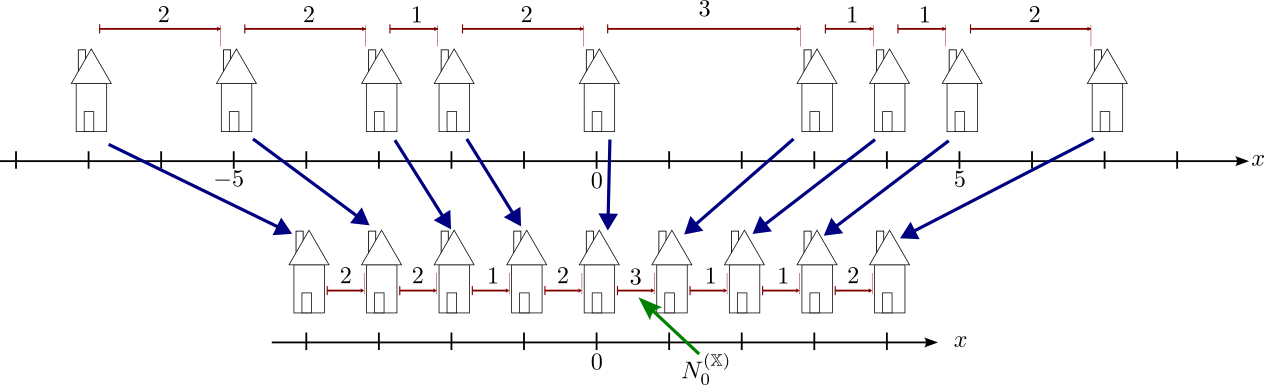}\caption{We fix the first existing house starting from $0$ as the new $x=0$
and line up all subsequent houses. Only the distance between houses
matters. The same can be analogously done for the lockdowns where
we only care about the total duration. \label{fig:alt-construction}}

\end{figure}

We use an alternative, more convenient description of the model. Instead
of considering Bernoulli random variables with parameters $q^{(t)}$
and $q^{(x)}$, we directly condense consecutive Bernoulli failures
into geometric random variables. Therefore, we will look at the total
duration of consecutive lockdowns instead of their existence at a
given time. Similarly, we consider distances between houses. The transition
from $\X_{x}$ to $\Nx_{x}$ is sketched in Figure \ref{fig:alt-construction}.
In terms of percolation, both constructions are equivalent. One just
loses information at which time step exactly a house recovers. 

\begin{defn}[Alternative construction]
\label{def:alternative-construction} Let $q^{(t)},q^{(x)},p\in(0,1)$
as well as $\alpha>1$ be given. We consider independent sequences
of independent geometric random variables $\Nt:=(\Nt_{t})_{t\in\Z}$
and $\Nx:=(\Nx_{x})_{x\in\Z}$ with parameters $q^{(t)}$ respectively
$q^{(x)}$.

Consider the graph $G=(\Z\times\Z,E)$ where $E$ consists of directed
edges of the form $(t,x)\to(t+1,y)$ with $t,x,y\in\Z$. We consider
a mixed bond-site-percolation model on $G$ where -- given $\Nt$
and $\Nx$ -- all vertices and edges are open (respectively closed)
independently from each other with probability
\begin{equation}
\P\{(t,x)\text{ is open}\,\vert\,\Nt\}=p^{\Nt_{t}}\label{eq:site-open}
\end{equation}
and
\[
\P\{(t,x)\to(t+1,y)\text{ is open}\,\vert\,\Nx\}=(1+d[x,y,\Nx])^{-\alpha}
\]
where $d[x,y,\Nx]$ is the distance between the $x$-th and $y$-th
house

\[
d[x,y,\Nx]:=\sum_{i=\min(x,y)}^{\max(x,y)-1}\Nx_{i}\,.
\]
\end{defn}

\begin{figure}

\includegraphics[width=1\columnwidth]{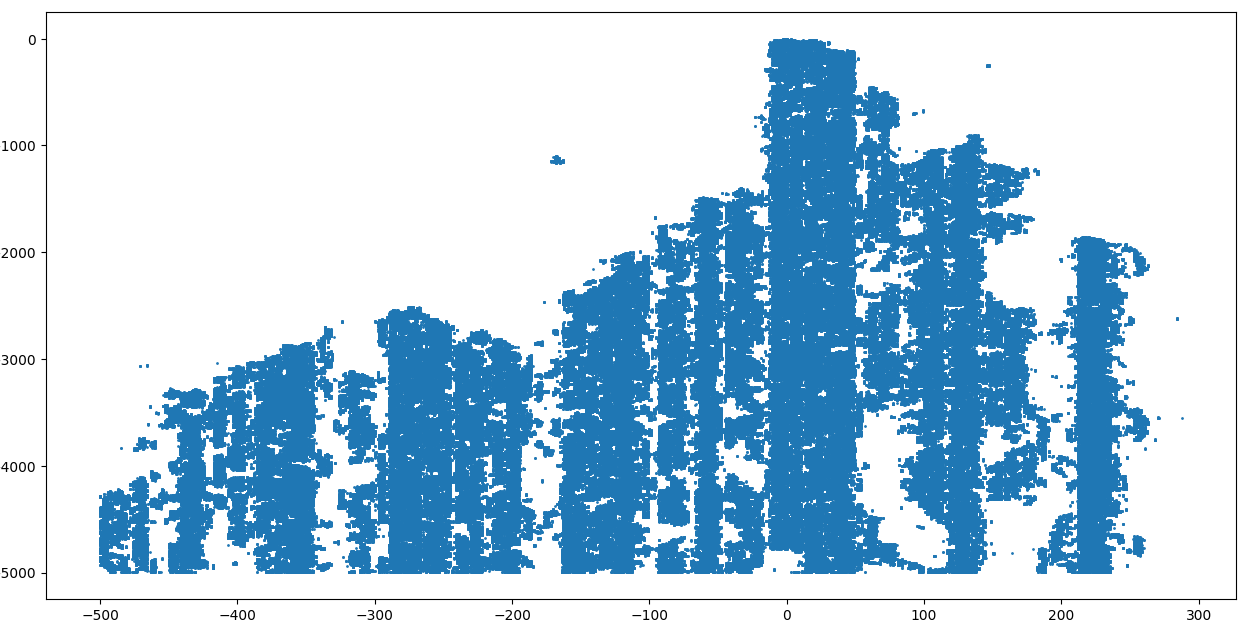}\caption{Simulation for $q^{(x)}=q^{(t)}=0.4$, $\alpha=3$ and $p=0.95$ starting
with an infected vertex in the origin. One can see distant infections
emerging due to long edges. Areas with large $\protect\Nx_{x}$ are
easily visible by the vertical gaps, but one can also see thin horizontal
gaps where $\protect\Nt_{t}$ is large. As the infection spreads in
space, it also seems to accelerate albeit often getting stuck at spatial
barriers. \label{fig:simulation}}
\end{figure}

One realisation of the condensed model is given in Figure \ref{fig:simulation}.
\begin{rem*}[Beyond geometric random variables]
 Note that for the alternative construction to make sense, we do
not actually need $\Nx,\Nt$ to be geometric random variables or even
to be $\N-$valued. In fact, it is perfectly reasonable to assume
$\Nx,\Nt\in\R_{>0}^{\Z}$ (which we will actually do in the following
rescaling lemmas).
\end{rem*}

The following two coupling lemmas allow us to freely choose the values
$q^{(t)}$ and $q^{(x)}$. We will be able to handle arbitrary $\alpha$
by choosing $p$ sufficiently large, so out of the four parameters
$q^{(t)},q^{(x)},\alpha,p$, we only need to focus on $p$.
\begin{lem}[Compensate $q^{(t)}$ by $p$]
\label{lem:Renormalize-qt} Let $\gamma>0$. Then, the $\MODEL$
with parameters $\gamma\Nt$ and $p^{1/\gamma}$ (with all other values
being unchanged) has the same distribution as the one with parameters
$\Nt,p$. In particular, we may assume that $q^{(t)}$ is arbitrarily
small by choosing $p$ accordingly close to $1$.
\end{lem}

\begin{proof}
This follows immediately from Equation (\ref{eq:site-open}).
\end{proof}
\begin{lem}[Compensate $q^{(x)}$ via $\alpha$]
\label{lem:Renormalize-qx} Let $\gamma\geq1$ and consider some
finite index set $J\subset\Z$. Then,
\[
\Big(1+\sum_{i\in J}\Nx_{i}\Big)^{-\alpha}\leq\Big(1+\sum_{i\in J}\lceil\gamma^{-1}\Nx_{i}\rceil\Big)^{-\gamma\alpha}\,.
\]
i.e., the $\MODEL$ with parameters $\Nx,\alpha$ is stochastically
dominated by the process with $\lceil\gamma^{-1}\Nx\rceil,\gamma\alpha$.
In particular, we may choose $q^{(t)}$ arbitrarily small by taking
$\alpha$ correspondingly large in order to show percolation.
\end{lem}

\begin{proof}
For every $a\geq0$, we prove $(1+a)^{\gamma}\geq1+\gamma a$. The
statement is true for $\gamma=1$. Differentiating in $\gamma$ at
$\gamma\geq1$ yields
\[
(1+a)^{\gamma}\cdot\log(1+a)\geq(1+a)\cdot a/(1+a)=a\,,
\]
so the statement holds for all $\gamma\geq1$. Finally,
\begin{align*}
\Big(1+\sum_{i\in J}\lceil\gamma^{-1}\Nx_{i}\rceil\Big)^{\gamma} & \geq1+\gamma\cdot\sum_{i\in J}\lceil\gamma^{-1}\Nx_{i}\rceil\geq1+\sum_{i\in J}\Nx_{i}
\end{align*}
which shows the claim after taking both sides to the power $-\alpha$.
\end{proof}

\subsection{Environment grouping scheme}

Next up is the grouping scheme for the random time and space environments.
Due to familiarity, we use the framework of \cite{MR2116736} rather
than \cite{MR4442896}. We extend it for more general values $\s,\cD$
and add extra details to the existing procedure. 

\medskip{}

\noindent\fbox{\begin{minipage}[t]{1\columnwidth - 2\fboxsep - 2\fboxrule}%
We fix two parameters
\begin{align*}
\s & \geq32\qquad\text{and}\qquad\cD<1/11\,.
\end{align*}
Consider stretches $N:=(N_{i})_{i\in\Z}$ with $N_{i}\in\N_{\geq1}\cup\{\infty\}$
where $N_{i}=\infty$ for at most one $i$.%
\end{minipage}}

\medskip{}

The bottom line is that, if the $N_{i}$ are generated by extremely
light-tailed iid geometric random variables, then the grouping scheme
terminates almost surely. As a reference, in \cite{MR2116736} we
have $\P(N_{i}\geq l+1)=(2^{-1000})^{l}$.
\begin{notation*}
From now on, $[m,n]$ will be an interval of integers, i.e.,
\begin{align*}
[m,n] & :=\{m,\,m+1,\dots,\,n-1,\,n\}\,,\\
(m,n) & :=[m,n]\backslash\{m,n\}\,.
\end{align*}
\end{notation*}
We group indices into bands depending on how ``bad'' they are. An
index $i\in\Z$ is bad if $N_{i}$ is large. These merge into bands
which are even ``worse''. We do so in a way such that bad bands
end up exponentially far apart. Unfortunately, a discount (depending
on the distance between far apart bands) has to be introduced for
the merging scheme to locally terminate almost surely for geometric
$N_{i}$.

\begin{figure}[th]
\includegraphics[width=1\columnwidth]{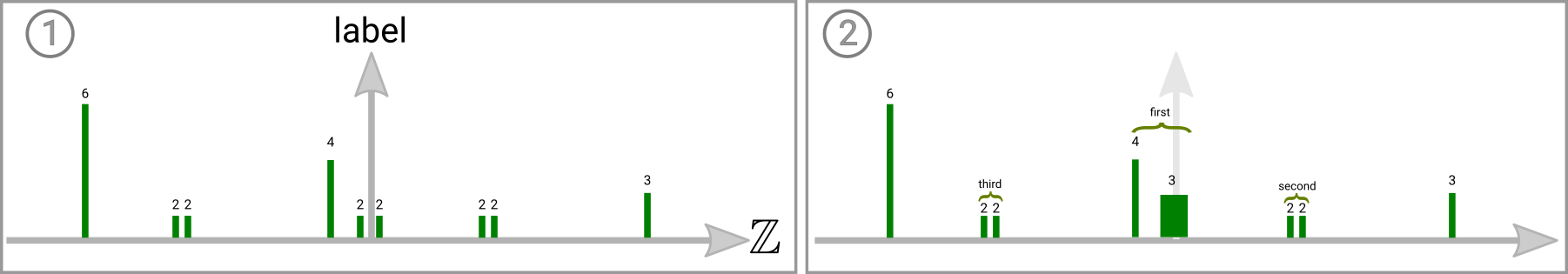}

\includegraphics[width=1\columnwidth]{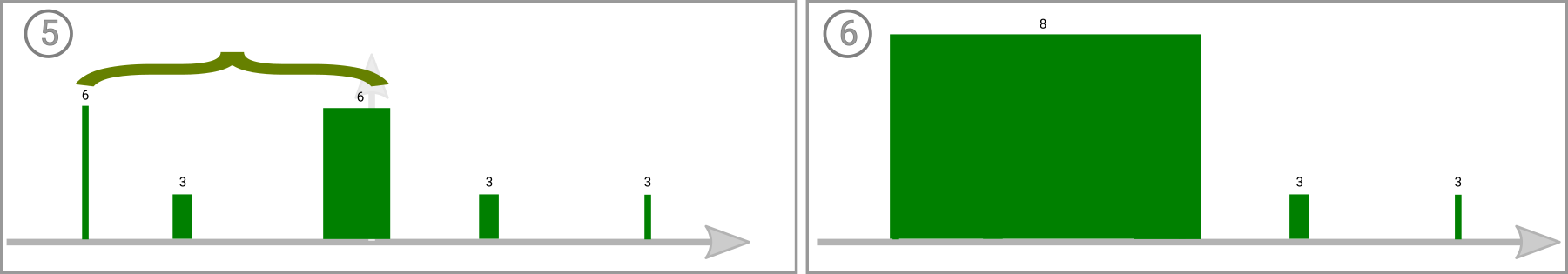}

\caption{\label{fig:-bands-and-labels}$k$ bands and labels for $k=1,2,5,6$.
The base height for labels in the diagrams is $1$. Curly brackets
show the merging order of the $k$ bands. After $k=6$, the merging
stops locally.}
\end{figure}

\medskip{}
We will consecutively define the $k$ bands of $N$, see Figure \ref{fig:-bands-and-labels}
for a rough illustration. 

\begin{defn}[$k$ bands and $k$ labels]
\label{def:bands-and-labels}\ The $k$ bands and $k$ labels are
defined inductively. A $1$ \textbf{band} is $\{i\}$ for $i\in\Z$.
The $1$ \textbf{label} of $\{i\}$ is
\[
f_{1}(i):=N_{i}\,.
\]
For indices $i,j\in\Z$, we set
\[
D_{k}(i,j):=\#\{k\text{ bands between }i\text{ and }j\text{ not containing either}\},
\]
e.g., at the current step $k=1$, we have $1+D_{1}(i,j)=|i-j|$. 

Given a partition of $\Z$ into $k$ bands together with their $k$
labels, the $k+1$ bands and $k+1$ labels are defined in the following
way: First, we pick specific \emph{merging indices} $i,j$ satisfying
\begin{equation}
\min\big(f_{k}(i),f_{k}(j)\big)-\log_{\s}\big(1+D_{k}(i,j)\big)>1\,.\label{eq:Distance-between-bands}
\end{equation}
The exact procedure for picking these is given in Algorithm \ref{alg:Merging-Indices}.
If no such pair exists, we terminate the merging scheme and set all
$k+1$ bands and labels to be the same as their $k$ counterpart.
Otherwise, using these $i,j$, we update as follows:
\end{defn}

\begin{enumerate}
\item Let $[m_{i},n_{i}]$ be the $k$ band containing $i$ and $[m_{j},n_{j}]$
the $k$ band containing $j$. Then, $[\tilde{m},\tilde{n}]$ is a
$k+1$ band with $\tilde{m}:=\min\{m_{i},m_{j}\}$ and $\tilde{n}:=\max\{n_{i},n_{j}\}$.
In this case, all $s\in[\tilde{m},\tilde{n}]$ have the $k+1$ \textbf{label}
\begin{equation}
f_{k+1}(s):=f_{k}(i)+f_{k}(j)-\left\lfloor \cD\log_{\s}\big(1+D_{k}(i,j)\big)\right\rfloor \,.\label{eq:label-update}
\end{equation}
Note that $f_{k+1}(s)\geq\max\{f_{k}(i),f_{k}(j)\}+2$.
\item Let $[\tilde{m},\tilde{n}]$ as above. If $[m,n]$ is a $k$ band
with $[m,n]\cap[\tilde{m},\tilde{n}]=\emptyset$, then it is also
a $k+1$ band. In this case, all $s\in[m,n]$ retain their label $f_{k+1}(s):=f_{k}(s)$.\\
Note that this condition is equivalent to $[m,n]\not\subset[\tilde{m},\tilde{n}]$.
\end{enumerate}
\begin{rem*}[Short summary]
 Each $k$ band is an interval of integers. At each step, two $k$
bands and everything inbetween merge into a bigger $k+1$ band of
larger label. In Algorithm \ref{alg:Merging-Indices}, we see that
$k$ bands close to the origin are preferred. For iid geometric $N_{i}$,
the merging procedure never terminates globally since there is always
something to merge.
\end{rem*}
Now, let us specify how exactly the merging indices in Definition
\ref{def:bands-and-labels} are chosen.
\begin{lyxalgorithm}[Finding merging indices]
\label{alg:Merging-Indices}\  Consider candidates $i,j\in\Z$ not
belonging to the same $k$ band and satisfying Equation \ref{eq:Distance-between-bands}.
\end{lyxalgorithm}

\begin{enumerate}
\item[1.] First, look for the smallest candidate pair $i,j$, that is, the
$i\in\Z$ with the smallest $|i+0.1|$ (i.e., $-|i|$ is preferred
over $|i|$) such that $|j|\leq|i|$. 
\item[2.] If $1+D_{k}(i,j)<(12\s)^{2}$, we choose $i,j$ as our merging indices.
\item[3.] If not, we try to look for ``better'' candidates that are close
to $i,j$: 
\begin{enumerate}
\item Search for candidates with $i',j'$ satisfying $1+D_{k}(i',j')<(12\s)^{2}$
as well as 
\[
1+D_{k}(i,j')<(12\s)^{2}\qquad\text{or}\qquad1+D_{k}(j,j')<(12\s)^{2}\,,
\]
 i.e.\ $j'$ is not too far away from $i$ or $j$, then continue
with $i',j'$ instead of $i,j$. (Note that $j'$ may coincide with
$i$ or $j$.)
\item If there are multiple candidates in the previous Step (a), take the
$j'$ minimizing $|j'+0.1|$ and then the $i'$ minimizing $D_{k}(i',j')$.
These are our merging indices.
\item If no such pair $i',j'\in\Z$ exists, take $i,j$ as the merging indices.
\end{enumerate}
\end{enumerate}
\begin{rem*}[Better candidates]
 The ``finding better candidates''-part is new compared to \cite{MR2116736}
and changes the order of merges. It is relevant for the proof of Theorem
\ref{lem:Main-Thm} Point 3 in the base case of simple bands (Definition
\ref{def:simple-bands}).
\end{rem*}
Two things are worth mentioning: First, if two $k$ bands with label
$\geq l$ are not at least $\s^{l-1}$ apart, then they will merge
at some point. Second, the size of a $k$ band (in terms of the indices
it contains) is limited by its label as seen in the following.
\begin{lem}[{Band size limit, \cite[Lemma 3.1]{MR2116736}}]
 \label{lem:Band-size-limit-by-label}If $[m,n]$ is a $k$ band
with $f_{k}(m)=l$, then $|n-m+1|\leq(\s/2)^{l-1}$.
\end{lem}

An indicated key result is the local termination of the merging scheme
for light-tailed $N_{i}$.

\begin{lem}[{Exponential decay of band labels, \cite[Lemma 3.4]{MR2116736}}]
\label{lem:Exp-Decay-Band-Labels} Assume that $N=(N_{i})_{i\in\Z}$
is a sequence of iid geometric random variables with $\P(N_{1}\geq l+1)=\mathfrak{q}^{l}$.
For any $J\in\Z$, $l\in\N$, and decay $\mathfrak{p}\in(0,1)$, there
exists a geometric parameter $\mathfrak{q}:=\mathfrak{q}(\s,\cD,\mathfrak{p})\in(0,1)$,
such that we have
\[
\P\big(\exists k\text{ s.t. }J\text{ lies in a }k\text{ band with label }\geq l\big)\leq\mathfrak{p}^{l-1}\,.
\]
In particular, the following holds almost surely: For each $J\in\Z$,
there exists a $K\in\N$ such that for all $k\geq K$, all the $n$
bands containing $J$ are identical.
\end{lem}

Since the $k$ bands are static at some point, we may now define the
``$k=\infty$'' bands.
\begin{defn}[Bands and labels]
\ 
\begin{enumerate}
\item An (integer) interval $[m,n]$ is called a \textbf{band} (without
$k$ in front) if there exists some $K\in\N$ such that $[m,n]$ is
a $k$ band for all $k\geq K$. For $j\in\Z$, the label of $j$ is
$f(j):=\lim_{k}f_{k}(j)$. The label of a band $[m,n]$ is $f(m)$.
\item If $N=(N_{i})_{i\in\Z}$ is such that $\Z$ decomposes into bands
that are finite, then we call $N$ \textbf{good}.
\end{enumerate}
\end{defn}

\medskip{}

Note that bands and their labels are always finite, i.e., $f(m)<\infty$,
except for the (potential) band containing $N_{i}=\infty$. 

From now on, we only deal with good $N=(N_{i})_{i\in\Z}$. 
\begin{cor}[]
\label{cor:Adding-Infinite-Label} In the setting of Lemma \ref{lem:Exp-Decay-Band-Labels},
we may with positive probability set $N_{0}=\infty$ without changing
the bands of $N$ and only changing the label of the band containing
$0$ to $\infty$.
\end{cor}

Setting $\Nt_{0}=\infty$ means that we consider all vertices of the
form $(0,x)$ to be closed. In this way, Corollary \ref{cor:Adding-Infinite-Label}
allows us to fix $0$ as a ``base height'' and therefore restrict
ourselves to a half space. 
\begin{defn}[Neighbouring bands and regularity]
\label{def:Enumeration-Neighbours}\ We enumerate bands as $B_{m}^{N},m\in\Z$
where $B_{0}^{N}$ is the band containing $0$ and $B_{1}^{N}$ is
the band to the right of $B_{0}^{N}$.
\begin{itemize}
\item Two bands $B_{m}^{N}$ and $B_{m'}^{N}$ are called \textbf{neighbouring
bands with labels} $\geq l$ if they both have labels $\geq l$ and
there is no band with label $\geq l$ inbetween.
\item The good sequence $N=(N_{i})_{i\in\Z}$ is called \textbf{regular}
if for all $l$ and all neighbouring bands $B_{m}^{N}$ and $B_{m'}^{N}$
with labels $\geq l$, we have $|m-m'|\in[\s^{l-1},\,12\cdot\s^{l-1})$,
i.e., there are at least $\s^{l-1}-1$ and at most $12\s^{l-1}-1$
bands between $B_{m}^{N}$ and $B_{m'}^{N}$.
\end{itemize}
\end{defn}

\medskip{}

A regular sequence is ``regular'' in the sense that bands with certain
labels show up regularly and are not spread too far apart. A good
sequence $N=(N_{i})_{i\in\Z}$ can always be made regular by artificially
raising individual $N_{i}$ (Lemma \ref{lem:Making-Sequences-Regular}).
We omit further details here since they are not needed to phrase the
general proof skeleton. The condition of $|m-m'|\geq\s^{l-1}$ is
automatically satisfied:
\begin{lem}[{\cite[Lemma 3.6]{MR2116736}}]
 If $B_{m}^{N}$ and $B_{m'}^{N}$ have label $\geq l$, $m\neq m'$,
then $|m-m'|\geq\s^{l-1}$.
\end{lem}

\begin{figure}
\includegraphics[width=1\columnwidth]{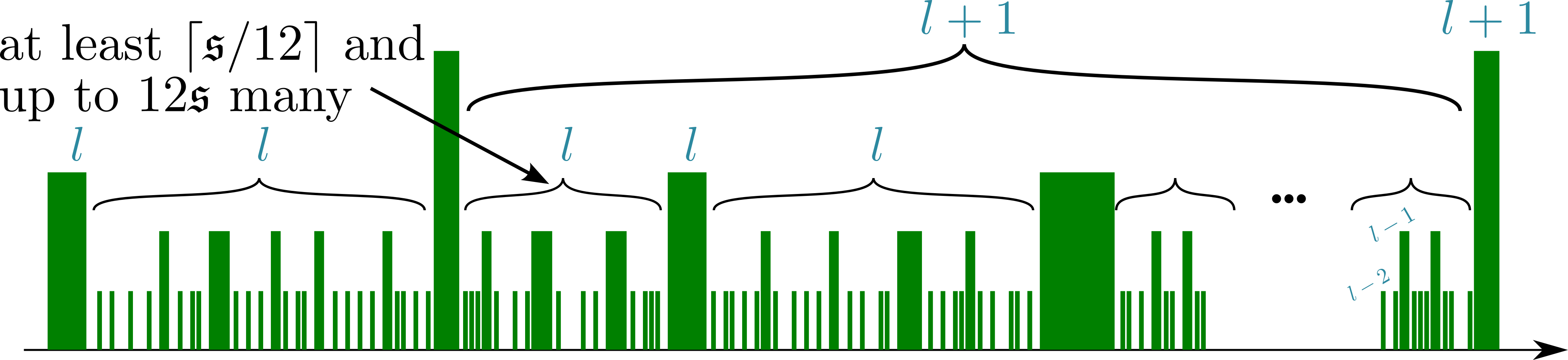}\caption{Bands (green bars) and segments (curly brackets) for regular $N$.
In this picture, there are always at least four $l$ segments between
two neighbouring bands of label $l$. \label{fig:regular}}
\end{figure}

\begin{proof}
If not, these bands would have merged before.
\end{proof}
Our next object of interest is ``the space between neighbouring bands''
since this is where our model will build up its ``bulk'' before
percolating through bands.
\begin{defn}[$l$ segments]
\label{def:Segment}\ Let $N$ be good and $[i_{1},\,i_{2}],\,[i_{3},\,i_{4}]$
be two neighbouring bands of label $\geq l$ (for $N$). Then we call
$(i_{2},i_{3})$ an $l$ \textbf{segment}. We refer to Figure \ref{fig:regular}
for an illustration of bands and segments for regular $N$.
\end{defn}

\begin{lem}[Number of $l$ segments between neighbouring bands]
\label{lem:Number-of-segments} Let $N$ be regular and $B_{m}^{N},\,B_{m'}^{N}$
be neighbouring bands of label $\geq l+1$. Let $\{m_{0},\dots,m_{k}\}=\{\tilde{m}\in[m,m']\,\vert\,B_{\tilde{m}}^{N}\text{ has label }\geq l\}$.
Then, $\lceil\s/12\rceil\leq k<12\cdot\s$. In particular, there are
between $\lceil\s/12\rceil$ and $12\cdot\s$ many $l$ segments separated
by bands of label $l$ between two neighbouring bands of label $\geq l+1$.
\end{lem}

\begin{proof}
Since $m_{i}-m_{i-1}<12\cdot\s^{l-1}$ by regularity and $m'-m=m_{k}-m_{0}\geq\s^{l}$,
we have $k\cdot12\cdot\s^{l-1}\geq\s^{l}$ which shows the first inequality.
The second follows from $m_{i}-m_{i-1}\geq\s^{l-1}$ and $m'-m=m_{k}-m_{0}<12\cdot\s^{l}$
by the same reasoning.
\end{proof}
Apart from the termination of the merging scheme (Lemma \ref{lem:Exp-Decay-Band-Labels}),
the above Lemma \ref{lem:Number-of-segments} is this section's important
take-away. It tells us that we always find a minimal amount of segments
between two bands. Regularity gives an upper bound.

\subsection{$n$ boxes in $\protect\Z\times\protect\Z$}

\begin{figure}
\includegraphics[width=1\columnwidth]{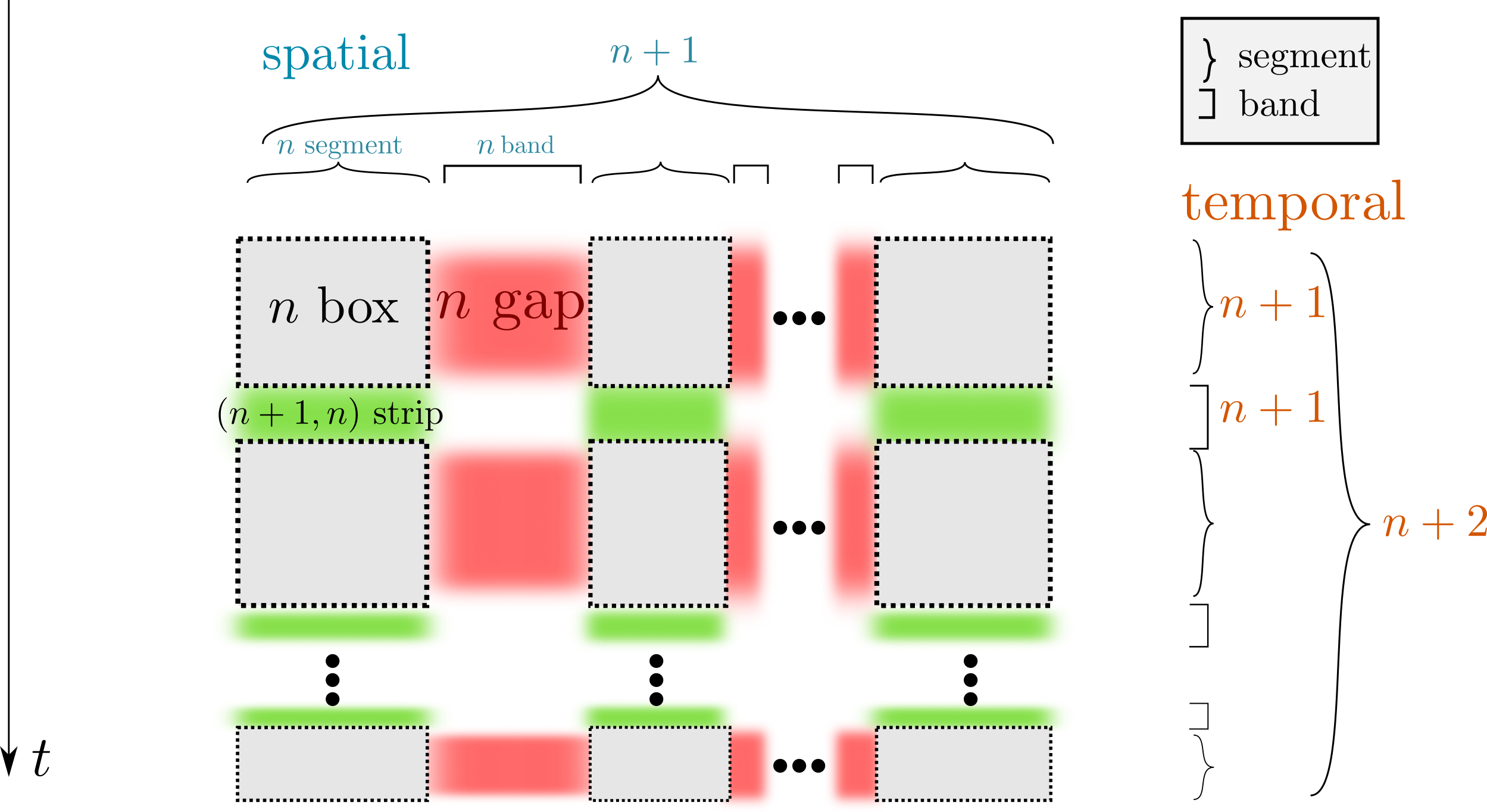}\caption{The inner structure of a $n+1$ box. We have grey $n$ boxes with
\textcolor{brown}{$n$ gaps} and \textcolor{teal}{$(n+1,n)$ strips}
between neighbouring boxes. Curly brackets depict segments while square
brackets depict bands. \label{fig:box-structure}}
\end{figure}

The framework for the environment grouping has been established. We
use it on the temporal environment with parameter $\s_{t}$ and the
spatial one with $\s_{x}$. Moving along our rough proof outline of
Section \ref{subsec:Idea-of-proof}, we now use this grouping to build
boxes. These boxes will be connected using ``inputs'' and ``outputs''
which are just vertices in special locations. 
\begin{defn}[$n$ boxes, $(m,n)$ strips and $n$ gaps]
 \ 
\begin{itemize}
\item If $[t_{1},t_{2}]$ is a temporal $2$ segment and if $\{x\}$ or
$\{x-1\}$ is a spatial band of label $1$, then any rectangle $[t_{1},t_{2}]\times\{x\}$
is a $1$ \textbf{box}. (Equivalently: if for every spatial band $[x_{1},x_{2}]$,
we have that $x\notin(x_{1},x_{2}]$.)
\item Let $n\in\N_{\geq2}$. Let $[t_{1},t_{2}]$ be a temporal $n+1$ segment,
i.e.\ the interval between two neighbouring bands with label $n+1$
(see Definition \ref{def:Segment}), and $(x_{1},x_{2})$ be a spatial
$n$ segment. Then, we call
\[
[t_{1},t_{2}]\times(x_{1},x_{2}]
\]
an $n$ \textbf{box}. (Yes, $x_{2}$ included!)
\item Let $n\in\N_{\geq1}$. Let $(t_{2},t'_{1})$ be a temporal $n+1$
band and $(x_{1},x_{2})$ be a spatial $m$ segment. Then, we call
\[
(t_{2},t'_{1})\times(x_{1},x_{2}]
\]
an $(n+1,m)$ \textbf{strip}. In other words: A $(n+1,n)$ strip is
the temporal interruption separating two vertically neighbouring $n$
boxes.
\item Let $n\in\N_{\geq1}$. Let $[t_{1},t_{2}]$ be a temporal $n+1$ segment
and $[x_{2},x_{1}']$ be a spatial $n$ band. Then, we call
\[
[t_{1},t_{2}]\times[x_{2},x_{1}']
\]
an $n$ \textbf{gap}. In other words: An $n$ gap is the spatial interruption
separating two horizontally neighouring $n$ boxes (starting at the
right-most border of the left box).
\end{itemize}
An illustration of an $n+1$ box is given in Figure \ref{fig:box-structure}.
\end{defn}

\begin{figure}
\includegraphics[width=1\columnwidth]{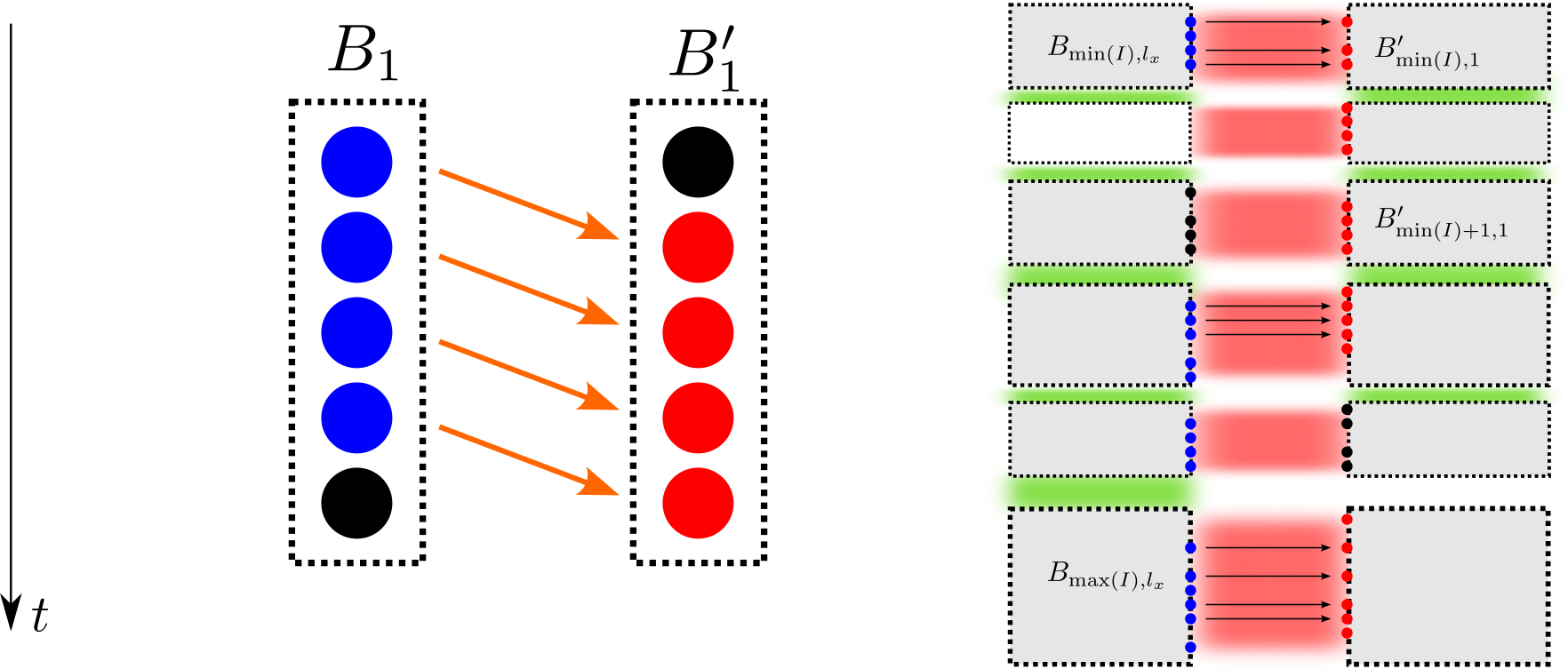}\caption{Connecting \textcolor{teal}{$\protect\Outh(B_{n+1})$} to \textcolor{brown}{$\protect\Inh(B_{n+1}')$}.
On the right, we can see that not all $n$ boxes are valid since their
respective $\protect\Inh/\protect\Outh$ might not be reachable. In
that picture, $B_{\min(I)+1,l_{x}}$ is bad (white). Therefore, that
sub-box as well as $B_{\min(I)+2,l_{x}}$ do not contribute to \textcolor{teal}{$\protect\Outh(B_{n+1})$}.
For \textcolor{brown}{$\protect\Inh(B_{n+1}')$}, it is only $B_{\max(I)-1,1}'$
that does not contribute since the gap to $B_{\max(I),1}'$ is bad,
i.e.\ $\protect\Inh(B_{\max(I)-1,1})$ might not be able to connect
to the rest of the $n+1$ box $B_{n+1}'$. \label{fig:io)}}
\end{figure}

\begin{rem*}[Renormalisation]
 We have to use $n+1$ rather than $n$ bands in the temporal part
because we essentially inserted a renormalisation step there. This
unfortunately also introduces a lot of bloat in notation. Lemma \ref{lem:Number-of-segments}
tells us that an $n+1$ box consists of between $\lceil\s_{x}/12\rceil+1$
and $12\s_{x}+1$ many columns as well as between $\lceil\s_{t}/12\rceil$
and $12\s_{t}$ many rows of $n$ boxes. These are separated by $n$
gaps respectively $(n,n)$ strips. \medskip{}
\end{rem*}
Next, we want to formally define good boxes as well as their inputs
and outputs now. The directed case makes things a bit more complicated,
but the multiscale arguments still work in a nice way. We will often
need to connect sets of vertices with each other, so it makes sense
to first introduce the following notion (slightly different to \cite{MR1901958}):
\begin{notation*}[(Fully) connected sets]
 Let $A,B\subset\Z^{2}$ be two sets of vertices. We write
\[
A\leadsto B
\]
 if there are $v\in A,w\in B$ such that $v\leadsto w$, i.e.\ there
exists an open directed path from $v$ to $w$. We write
\[
A\fc B
\]
if for every $v\in A$ and every $w\in B$, we have $v\leadsto w$.
Note that 
\[
A\fc B\leadsto C\fc D\implies A\fc D\,.
\]
\end{notation*}
\begin{rem*}
Before directly moving on to the definition of inputs and outputs,
let us recall the basic idea first. Each ``good'' $n$ box $B_{n}$
will have four sets of vertices $\Inv(B_{n})$ (on the top), $\Inh(B_{n})$
(on the sides), $\Outh(B_{n})$ (also on the sides) and $\Outv(B_{n})$
(on the bottom). The $\Out$ stands for outgoing connections to other
boxes' ingoing connections $\In$. For example, $\Outv(B_{n})$ stands
for vertices which potentially build an open path to $\Inv(B_{n}')$
for another $n$ box $B_{n}'$ directly below $B_{n}$. Since the
cardinality of these sets grows exponentially in $n$, this means
we will have exponentially many trials to bridge an $(n,n)$ strip
(and analogously $n$ gaps).
\end{rem*}
The locations of $\Inh(B_{n})$ and $\Outh(B_{n})$ have to be set
carefully so that the inputs and outputs are sufficiently well connected
inside $B_{n}$. Furthermore, we are only able to make statements
on ``good'' $n$ boxes, so the following definition will appear
quite bloated.
\begin{defn}[Good $n$ boxes, inputs and outputs]
\label{def:Good-boxes} Let $B_{n}$ be an $n$ box.
\begin{itemize}
\item For $n=1$, the $n$ box $B_{n}=[t_{1},t_{2}]\times\{x\}$ is \textbf{good}
if all vertices are open (in the sense of Definition \ref{def:alternative-construction}).
In this case, we write $\Inv(B_{n}):=\{(t_{1},x)\},\Outv(B_{n}):=\{(t_{2},x)\}$
as well as
\[
\Inh(B_{n}):=(t_{1},t_{2}]\times\{x\}\qquad\text{and}\qquad\Outh(B_{n}):=[t_{1},t_{2})\times\{x\}\,.
\]
\item An $n$ \textbf{gap} between two horizontally neighbouring boxes $B_{n},B_{n}'$
is \textbf{good} if 
\[
\Outh(B_{n})\leadsto\Inh(B_{n}')\qquad\text{and}\qquad\Outh(B_{n}')\leadsto\Inh(B_{n})\,.
\]
\item An $(n,n)$ \textbf{strip} between two vertically neighbouring boxes
$B_{n},B_{n}'$ is \textbf{good} if
\[
\Outv(B_{n})\leadsto\Inv(B_{n}')\,.
\]
\item We call an $n+1$ box $B_{n+1}$ \textbf{good} (and otherwise \textbf{bad})
if the sum of the number of the following bad objects is at most $1$:
\begin{itemize}
\item[A)]  $n$ boxes inside $B_{n+1}$,
\item[B)]  $(n+1,n)$ strips between two $n$ boxes inside $B_{n+1}$,
\item[C)]  $n$ gaps between two $n$ boxes inside $B_{n+1}$.
\end{itemize}
\item In the case of $B_{n+1}$ being good, we first number its $n$ boxes
$(B_{i,j})_{1\leq i\leq l_{t},1\leq j\leq l_{x}}$ by their location
(with $i=j=1$ being top-left) where $l_{t}\in[\lceil\s_{t}/12\rceil,\,12\s_{t}]$
and analogously $l_{x}-1\in[\lceil\s_{x}/12\rceil,\,12\s_{x}]$. Next,
we set (for some $\khor\in\N$ specified in Equation (\ref{eq:kappa_ver-and-kappa_hor}))
\[
I:=[0,\khor+4)+12\s_{x}+2\,.
\]
Then, we can finally define the \textbf{inputs} and \textbf{outputs}
of the $n+1$ box. The vertical inputs/outputs are as follows: For
$j\in\{1,\dots,l_{x}\}$, we set
\begin{align*}
\Inv(B_{n+1}) & :=\left\{ v\in\Inv(B_{1,j})\,\vert\,B_{1,j}\text{ is good}\right\} \\
\Outv(B_{n+1}) & :=\left\{ v\in\Outv(B_{l_{t},j})\,\vert\,B_{l_{t},j}\text{ is good }\right\} \,.
\end{align*}
Let $\partial B_{n+1}\subset\Z\times\Z$ be the boundary, i.e.\ the
set of all vertices in $B_{n+1}$ having a neighbour outside of it.
Then,
\begin{align*}
\Inh(B_{n+1}) & :=\big\{ v\in\Inh(B_{i,j})\cap\partial B_{n+1}\,\vert\,B_{i,j},B_{i+1,j}\text{ are valid, }j\in\{1,l_{x}\},i\in I\big\}\\
\Outh(B_{n+1}) & :=\big\{ v\in\Outh(B_{i,j})\cap\partial B_{n+1}\,\vert\,B_{i-1,j}B_{i,j}\text{ are valid, }j\in\{1,l_{x}\},i\in I\big\}
\end{align*}
where we say that $n$ boxes $B_{n},B_{n}'$ are \textbf{valid} if
both are good and $\Outv(B_{n})\leadsto\Inv(B_{n}')$. 
\end{itemize}
We refer to Figure \ref{fig:io)} for an illustration.
\end{defn}

The parameters $\s_{x}$ and $\s_{t}$ roughly correspond to the width
respectively height of the given boxes. Thus, they also influence
the number of connectors between boxes: The larger $\s_{x}$, the
larger the number of vertical connectors between vertically neighboured
boxes (since the boxes are wider). The same holds for $\s_{t}$. We
will capture the minimal amount of vertical (respectively horizontal)
connectors via the parameters $\kver$ and $\khor$.

\medskip{}

\noindent\fbox{\begin{minipage}[t]{1\columnwidth - 2\fboxsep - 2\fboxrule}%
We set the following parameters:
\begin{equation}
\kver:=\lceil\s_{x}/12\rceil-2\qquad\text{and}\qquad\khor:=\lceil\s_{t}/12\rceil-2\cdot(12\s_{x}+1)-4\,\label{eq:kappa_ver-and-kappa_hor}
\end{equation}
and assume $\kver\geq64$ (additional conditions on $\khor$ are specified
later).%
\end{minipage}}

\medskip{}

\begin{rem*}
The spatial parameter $\s_{x}$ can just be fixed to $12\cdot66$
to ensure $\kver\geq64$. The value of $\s_{t}$ (equivalently $\khor$)
will however depend $\alpha$ and a small parameter $\pp$ which governs
the probability of bad boxes introduced later in Lemma \ref{lem:Main-Thm}.
Also, for a rough estimate on the values: We already have 
\begin{align*}
\s_{x} & \geq11\cdot66+1\geq700\qquad\text{and}\qquad\s_{t}\geq17'000\,,
\end{align*}
so this is quantitively unfeasible.
\end{rem*}

\subsection{Towards proving percolation}

The parameters $\kver,\khor$ had to be set in such a convoluted way
to ensure the following connectivity inside good $n$ boxes:
\begin{lem}[Connecting inputs and outputs inside]
\label{lem:Connecting-Inputs-Outputs-Inside} Let $n\in\N$. Let
$B_{n}$ be a good $n$ box. Then,
\begin{align*}
\Inv(B_{n}) & \fc\Outh(B_{n})\\
\Inv(B_{n}) & \fc\Outv(B_{n})\\
\Inh(B_{n}) & \fc\Outv(B_{n})\,.
\end{align*}
In particular, if $B_{n}'$ is a horizontally neighbouring good $n$
box with the $n$ gap inbetween being good as well, then 
\[
\Inv(B_{n})\fc\Outv(B_{n}')\,.
\]
\end{lem}

The application of this lemma can be retrospectively seen in Figure
\ref{fig:strategy}.

We may finally state the main auxilliary theorem for the multiscale
argument. Using the probability of good $n$ boxes, we are then in
the state to prove the main theorem on the survival of the infection
(Theorem \ref{thm:Survival}).
\begin{lem}[{Main auxilliary lemma, \cite[Lemma 4.3]{MR2116736}}]
\label{lem:Main-Thm}\ Let $\pp\in(0,1)$ and $\kver\geq64$. For
all sufficiently large $\khor\in\N$ (depending on $\pp,\kver$),
there exists $p_{c}\in(0,1)$ such that in the $\MODEL$ model for
any $p\geq p_{c}$
\begin{enumerate}
\item $\P(B_{n}\text{ is good})\geq1-\pp^{n+1}$ for any $n$ box $B_{n}$.
\item Let $G_{n}$ be a temporal $n$ gap (between two neighbouring $n$
boxes). Then,
\[
\P(G_{n}\text{ is good})\geq1-\pp^{n+1}\,.
\]
\item For an $(n+1,n)$ strip $\bar{S}$ between two $n$ boxes $B_{n},B_{n}'$,
we have
\[
\P(\bar{S}\text{ is good})\geq1-\pp^{n+1}\,.
\]
\end{enumerate}
\end{lem}

\begin{proof}[Proof outline]
 For the reader's convenience, we will give a brief overview over
the main steps. The complete proof will be given in Section \ref{sec:Proving-Percolation}.
\begin{itemize}
\item Point 1 follows from combinatorial estimates (Lemma \ref{lem:Multiscale-Estimate})
and induction after proving Point 2 and 3. 
\item $n$ gaps are exponentially large in $n$ (Lemma \ref{lem:maximal-stretch-of-simple}).
Since we have long-range edges, we can guarantee Point 2 by choosing
$\khor$ large depending on $\alpha$ (Lemma \ref{lem:Gap-Crossing}).
We do so by crossing the whole $n$ gap in a single jump. 
\item Point 3 follows from the main difficulty of the whole procedure: the
``drilling'' (Proposition \ref{prop:Drilling}). Luckily, the proof
in \cite{MR2116736} still works here.
\end{itemize}
\end{proof}
Taking Lemma \ref{lem:Main-Thm} Point 1, we can finally prove the
existence of an infinite directed path and in particular the phase
transition of the $\MODEL$ in the parameter $p$.
\begin{proof}[Proof of Theorem \ref{thm:Survival}]
 \ Puzzling everything together is still something.
\begin{enumerate}
\item We first take $\s_{x}=12\cdot66,\,\pp=1/4$.
\item Using Lemma \ref{lem:Renormalize-qx}, we may assume at the cost of
$\alpha$ that $q^{(x)}$ is sufficiently small such that Lemma \ref{lem:Exp-Decay-Band-Labels}
holds for $\cD=1/12$, in particular we may use the whole framework
of Section \ref{sec:Bands-Labels-Regularity}. 
\item Lemma\ \ref{lem:Main-Thm} gives us some $\s_{t}$ and $p_{c}$ for
which it holds.
\item Using Lemma \ref{lem:Renormalize-qt}, we may assume at the cost of
$p$ that $q^{(t)}$ is sufficiently small such that Section \ref{sec:Bands-Labels-Regularity}
can be used for that $\s_{t}$ and $\cD$. 
\item Corollary \ref{cor:Adding-Infinite-Label} lets us fix base height
$0$ for a positive fraction of temporal environments, i.e.\ $\Nt_{0}=\infty$.
\item Next, choose $u=(1,42)\in\Z_{>0}\times\Z$. This lies in some $n$
box for $n$ large enough. By Lemma\ \ref{lem:Main-Thm} Point 1
and Borel--Cantelli, there exists some $N_{0}$ such that all the
$n$ boxes $B_{n}$ with $n\geq N_{0}$ containing $u$ are good. 
\item Now, take any $v\in\Inv(B_{N_{0}})$. Since $\Nt_{0}=\infty$, we
have $v\in\Inv(B_{n})$ for every $n\geq N_{0}$, in particular $v\fc\Out(B_{n})$.
Therefore, $v\leadsto w$ for infinitely many $w$. This already yields
us an infinite directed path: Set $v_{0}:=v$. Since $v_{0}$ only
has finitely many direct successors, we may choose any of these successors
$v_{1}$ that has infinitely many $w$ with $v_{0}\to v_{1}\leadsto w$.
Inductively continuing this scheme, we obtain an infinite path $v_{0}\to v_{1}\to v_{2}\to\dots$.
\item $\Nt=(\Nt_{i})_{i\in\Z}$ is an iid sequence, in particular ergodic.
So
\[
\P\{\Nt\,\text{s.t.}\,\P(\exists v\in\Z^{2},\,v\leadsto\infty\,\vert\Nt)=1\}\in\{0,1\}\,.
\]
Since we have proven percolation on a positive fraction of environments,
it has to hold for almost all of them.
\end{enumerate}
\end{proof}

\section{Details: environment grouping \label{sec:Bands-Labels-Regularity}}

Now, that the rigorous roadmap has been laid out in Section \ref{sec:Proof-skeleton},
it is time to flesh it out. The main goals in the current sections
are:
\begin{itemize}
\item Showing Lemma \ref{lem:Exp-Decay-Band-Labels}, i.e., local termination
of the merging scheme.
\item Showing how good sequences can always be made regular and even ``very
regular'' (Lemmas \ref{lem:Making-Sequences-Regular}, \ref{lem:Regular-to-Very-Regular}).
\item Introducing the notion of simple bands (Definition \ref{def:simple-bands})
and how they are well-behaved (Lemma \ref{lem:maximal-stretch-of-simple}).
\item Splitting up very regular bands (Lemma \ref{lem:q-estimates}).
\end{itemize}

\subsection{Local termination of merging scheme}

We start by quantifying the maximal ``size'' of bands, i.e., giving
the proof for Lemma \ref{lem:Band-size-limit-by-label}. 
\begin{proof}[Proof of Lemma \ref{lem:Band-size-limit-by-label}]
The statement is true for $l\leq3$ since it implies $m=n$. Now
suppose $[m,n]$ is a $k$ band with $m\neq n$ and $f_{k}(j)=l>3$.
Then, there must be some $k'<k$ and $m',n'$ such that the $k'$
bands $[m,m']$ and $[n',n]$ merge into $[m,n]$. We denote $\underline{l}:=f_{k'-1}(m),\,\bar{l}=f_{k'-1}(n)$
and $\mathbf{N}:=D_{k'}(m',n')$. Then, there are at most $\mathbf{N}/\s^{L-1}$
many $k'$ bands with labels $L$ between $[m,m']$ and $[n',n]$
(otherwise some would have merged). Using the induction hypothesis
on the $k'$ bands of labels $L$
\begin{align*}
|n- & m+1|\leq|m'-m+1|+|n'-m'-1|+|n-n'+1|\\
\leq & (\tfrac{\s}{2})^{\underline{l}-1}+(\tfrac{\s}{2})^{\bar{l}-1}+\sum_{L=1}\sum\big\{|b'-b+1|\,\Big|\,[b',b]\subset(m',n')\text{ is }k'-1\text{ band with label }L\big\}\\
\leq & (\tfrac{\s}{2})^{\underline{l}-1}+(\tfrac{\s}{2})^{\bar{l}-1}+\sum_{L}(\mathbf{N}/\s^{L-1})\cdot(\tfrac{\s}{2})^{L-1}\leq(\tfrac{\s}{2})^{\underline{l}-1}+(\tfrac{\s}{2})^{\bar{l}-1}+2\mathbf{N}\\
^{(*)}\leq & (\tfrac{\s}{2})^{\max\{\underline{l},\bar{l}\}-1}+(\tfrac{\s}{2})^{\min\{\underline{l},\bar{l}\}-1}+2\cdot\s^{\min\{\underline{l},\bar{l}\}-1}\leq\tfrac{4}{\s}\cdot(\tfrac{\s}{2})^{\max\{\underline{l},\bar{l}\}}+\tfrac{2}{\s}\cdot\s^{\min\{\underline{l},\bar{l}\}}\\
^{(**)}\leq & \tfrac{4}{\s}(\tfrac{\s}{2})^{l-1}+\tfrac{2}{\s}\cdot\s^{\min\{\underline{l},\bar{l}\}}\leq\tfrac{4}{\s}(\tfrac{\s}{2})^{l-1}+\tfrac{2}{\s}\cdot\s^{l/(2-\cD)}\\
^{(***)}\leq & \tfrac{4}{\s}(\tfrac{\s}{2})^{l-1}+\tfrac{1}{2}\cdot(\tfrac{\s}{2})^{l-1}\leq(\tfrac{\s}{2})^{l-1}\,,
\end{align*}
where $(*)$ follows from Equation \ref{eq:Distance-between-bands}
being equivalent to $N<\s^{\min\{\underline{l},\bar{l}\}-1}-1$, $(**)$
uses $l\geq\max\{\underline{l},\bar{l}\}+1$ and $(***)$ uses Equation
\ref{eq:label-update} for
\begin{align*}
l & =\max\{\underline{l},\bar{l}\}+\min\{\underline{l},\bar{l}\}-\cD\log_{\s}(1+\mathbf{N})\\
 & \geq2\min\{\underline{l},\bar{l}\}-\cD(\min\{\underline{l},\bar{l}\}-1)\geq(2-\cD)\min\{\underline{l},\bar{l}\},
\end{align*}
which yields
\[
4\cdot\s^{l/(2-\cD)}\leq(\tfrac{\s}{2})^{l}\,,
\]
(together with $\cD<1/11\leq2-(1-\log_{32}4)^{-1}\leq2-(1-\log_{\s}4)^{-1}$).
\end{proof}

\begin{cor}[{Combining distant bands, \cite[Lemma 3.2]{MR2116736}}]
\label{cor:Combine-distant-bands} If $[m,m']$ and $[n',n]$ merge
at step $k+1$, then
\[
\mathbf{N}:=\#\big\{ k\text{ bands lying in }(m',n')\big\}\geq\tfrac{2}{3}(n'-m'-1)\,.
\]
\end{cor}

\begin{proof}
Let $\mathbf{N}_{l}$ be the number of $k$ bands in $(m',n')$ with
label $l$. Then $\mathbf{N}_{l}\leq\mathbf{N}/(\s^{l-1})$ since
otherwise some would have merged first. Furthermore $\mathbf{N}=\sum_{l}\mathbf{N}_{l}$,
so 

\begin{align*}
n'-m'-1 & =\sum\big\{ b'-b+1\,\vert\,[b,b']\text{ is a }k\text{ band in }(m',n')\big\}\\
\leq\mathbf{N}_{1}+ & \mathbf{N}_{2}+\sum_{l\geq3}\mathbf{N}/(\s^{l-1})(\s/2)^{l-1}\leq\mathbf{N}_{1}+\mathbf{N}_{2}+\mathbf{N}/2\leq\tfrac{3}{2}\mathbf{N,}
\end{align*}
which shows the claim.
\end{proof}
Next we consider generators. They are relevant in this subsection
as well as in Section \ref{subsec:Regular-Bands}. Generators of a
$k$ band are, loosely speaking, the boundary points of $<k$ bands
that are merged to form the final $k$ band.
\begin{defn}[Generators of a $k$ band]
 Let $[m,n]$ be a $k$ band.
\begin{itemize}
\item The $k$ \textbf{generators} of $[m,n]$ are $m$ and $n$.
\item For $k'<k$, the $k'$ \textbf{generators} of $[m,n]$ are the $k'$
generators of the $k'$ bands containing a $k'+1$ generator of $[m,n]$.
\item For $1$ generators, we will omit the $1$ and just call them \textbf{generators}.
\item We call a generator $g$ a \textbf{maximal generator} if it satisfies
the following:
\begin{itemize}
\item[] If the $k'<k$ bands $[m_{1},n_{1}]$, $[m_{2},n_{2}]$ with $g\in[m_{1},n_{1}]$
combine, then $f_{k'}(m_{1})\geq f_{k'}(m_{2})$.
\end{itemize}
\item One verifies that $k$ bands $B$ always have a maximal generator.
Pick the smallest one $\cg_{k}(B)$.
\end{itemize}
The next lemma limits the possibilities of generators to be spread
apart.
\end{defn}

\begin{lem}[{\cite[Lemma 3.3]{MR2116736}}]
\label{lem:log-dist-of-generators} Let $i_{1}<i_{2}<\dots<i_{n}$
be the generators of an $k$ band with 
\[
\sum_{j=1}^{n}f_{1}(i_{j})=m\,.
\]
Then, there exists $C(\s,\cD)>0$ such that
\[
\sum_{j=2}^{n}\lfloor\log_{2}(i_{j}-i_{j-1}+1)\rfloor\leq C(\s,\cD)m\,.
\]
\end{lem}

\begin{proof}
If $n=1$, we are done. For each $j\in[2,n]$, let $k_{j}$ be the
value such that there exists $a$ and \textbf{$b$} such that the
$k_{j}$ bands $[i_{a},i_{j-1}]$ and $[i_{j},i_{b}]$ merge into
$[i_{a},i_{b}]$. Let $\mathbf{N}_{j}$ be number of $k_{j}$ bands
between $[i_{a},i_{j-1}]$ and $[i_{j},i_{b}]$. By the previous Corollary
\ref{cor:Combine-distant-bands}
\[
\tfrac{2}{3}(i_{j}-i_{j-1}+1)\leq1+\mathbf{N}_{j}.
\]
 We have that $n\leq m/2$ as well as
\begin{equation}
\sum_{j=2}^{n}\lfloor\cD\log_{\s}(1+\mathbf{N}_{j})\rfloor\leq m\,,\label{eq:idk1}
\end{equation}
whose proof will be given right after. $\lfloor xy\rfloor\geq x\lfloor y\rfloor-1$
yields the following chain of implication
\[
m\geq\sum_{j=2}^{n}\lfloor\cD\log_{\s}(\tfrac{2}{3}(i_{j}-i_{j-1}+1))\rfloor\geq\sum_{j=2}^{n}\cD\log_{\s}2\cdot\lfloor\log_{2}(\tfrac{2}{3}(i_{j}-i_{j-1}+1))\rfloor-n
\]
\[
m\frac{3\log_{2}\s}{2\cD}\geq\sum_{j=2}^{n}\lfloor\log_{2}(\tfrac{2}{3}(i_{j}-i_{j-1}+1))\rfloor\geq\sum_{j=2}^{n}\lfloor\log_{2}(i_{j}-i_{j-1}+1)\rfloor-n
\]
\[
m\big(\frac{3\log_{2}\s}{2\cD}+\frac{1}{2}\big)\geq\sum_{j=2}^{n}\lfloor\log_{2}(i_{j}-i_{j-1}+1)\rfloor\,.
\]
The claim follows from choosing $C(\s,\cD):=\big(\frac{3\log_{2}\s}{2\cD}+\frac{1}{2}\big)$.
\end{proof}
\begin{proof}[Proof of Inequality (\ref{eq:idk1})]
 By the assumption $\sum_{j=1}^{n}f_{1}(i_{j})=m$, we first see
that $n\leq m/2$ since $f_{1}(i_{j'})\geq2$ for generators. Equation
(\ref{eq:label-update}) gives
\[
f_{k_{j}+1}(i_{j})=f_{k_{j}}(i_{j-1})+f_{k_{j}}(i_{j})-\lfloor\cD\log_{\s}(1+\mathbf{N}_{j})\rfloor\,.
\]
If for example $\min_{j}k_{j}=k_{j'}$ is the smallest, i.e., we first
combine $[i_{j'-1}]$ and $[i_{j'}]$, this would yield
\begin{align*}
m & =\sum_{j=1}^{n}f_{1}(i_{j})=\sum_{j=1}^{n}f_{k_{j'}}(i_{j})=\sum_{j=2}^{n}f_{k_{j'}+1}(i_{j'})+\lfloor\cD\log_{\s}(1+\mathbf{N}_{j'})\rfloor\,.
\end{align*}
Continuing iteratively with $f_{k_{j'}+1}$ now instead of $f_{1}$
yields
\begin{align*}
m & =f_{k'}(i_{j'})+\sum_{j=2}^{n}\lfloor\cD\log_{\s}(1+\mathbf{N}_{j})\rfloor\geq\sum_{j=2}^{n}\lfloor\cD\log_{\s}(1+\mathbf{N}_{j})\rfloor\,,
\end{align*}
which finishes the calculation. (Even $f_{k'}(i_{j'})\geq2n$ since
merges raise labels by at least $2$.)
\end{proof}
We are finally in the spot to prove the first milestone: local termination
of the merging scheme.

\begin{proof}[Proof of Lemma \ref{lem:Exp-Decay-Band-Labels}]
 Let $\mathfrak{q}^{1/2}\leq\mathfrak{p}\cdot\s^{-3C(\s,\cD)}/2$
and $l>1$. Assume that $J$ actually lies in a $k$ band with label
$\geq$l and generators $i_{1}<\dots<i_{n}$. In the case that $i_{1}=i_{n}=J$,
the claim follows
\[
\P(f_{1}(J)\geq l)=\P(N_{J}\geq l)=\mathfrak{q}^{l-1}\leq\mathfrak{p}^{l-1}\,.
\]
Otherwise, we continue. The $i_{j}$ satisfy by the label updating
procedure in Definition \ref{def:bands-and-labels}
\[
m:=\sum_{j}^{n}f_{1}(i_{j})\geq l\,.
\]
By Lemma \ref{lem:log-dist-of-generators} above and Lemma \ref{lem:Combinations-of-sums}
below, we have at most $2^{C(\s,\cD)m}$ choices for $\lfloor\log_{2}(i_{2}-i_{1}+1)\rfloor,\dots,\lfloor\log_{2}(i_{n}-i_{n-1}+1)\rfloor$.

Given one such choice, we yet again have $2^{(C(\s,\cD)+\tfrac{1}{2})m}$
choices for $(i_{2}-i_{1}),\dots(i_{n}-i_{n-1})$:

\medskip{}

\noindent\fbox{\begin{minipage}[t]{1\columnwidth - 2\fboxsep - 2\fboxrule}%
Set $a_{j}:=\lfloor\log_{2}(i_{j}-i_{j-1}+1)\rfloor$, in particular
$i_{j}-i_{j-1}\leq2^{a_{j}+1}$. There are $2^{a_{j}+1}$ possibilities
for each individual $j$, so in total for the whole ensemble
\[
\prod_{j=2}^{n}2^{a_{j}+1}\leq2^{C(\s,\cD)m+n}\leq2^{\{C(\s,\cD)+1/2\}m}\,.
\]
\end{minipage}}

\medskip{}
Furthermore, there are at most $(\s/2)^{l-1}$ possible starting
locations for $i_{1}$ since by Lemma \ref{lem:Band-size-limit-by-label}
\[
i_{1}\leq J\leq i_{1}+(\s/2)^{l-1}-1\,.
\]
So in total, we have at most $(\s/2)^{l-1}\cdot2^{(2C(\s,\cD)+1)m}$
choices for $i_{1},\dots i_{n}$. For each choice of $i_{1},\dots,i_{n}$,
there are at most $2^{m}$ choices for $f_{1}(i_{1}),\dots f_{1}(i_{n})$
(by Lemma \ref{lem:Combinations-of-sums} below), so we have at most
\[
(\s/2)^{l-1}\cdot2^{(2C(\s,\cD)+1)m}\cdot2^{m}\leq\s^{3C(\s,\cD)m}
\]
choices for the combined $i_{j}$ and $f_{1}(i_{j})$. Each such choice
has probability $\mathfrak{q}^{m-n}\leq\mathfrak{q}^{m/2}$ since
$\P\{f_{1}(i_{j})=s\}\leq\mathfrak{q}^{s-1}$. Therefore, for $\mathfrak{q}^{1/2}\leq\mathfrak{p}\cdot\s^{-3C(\s,\cD)}/2$
(in particular $\mathfrak{q}\leq\mathfrak{p}/2$)
\begin{align*}
\P\big( & \exists k\text{ s.t. }j\text{ lies in an }k\text{ band with label }\geq l\big)\leq\sum_{m\geq l}\big[\s^{3C(\s,\cD)}\mathfrak{q}^{1/2}\big]^{m}+\mathfrak{q}^{l-1}\\
 & \leq\sum_{m\geq l}(\mathfrak{p}/2)^{m}+(\mathfrak{p}/2)^{l-1}=\mathfrak{p}^{l}\frac{1}{2^{l}(1-\mathfrak{p}/2)}+(\mathfrak{p}/2)^{l-1}\leq2^{-(l-1)}\cdot\big[\mathfrak{p}^{l}+\mathfrak{p}^{l-1}\big]\leq\mathfrak{p}^{l-1},
\end{align*}
as desired.
\end{proof}
Here is the auxiliary lemma we previously used.
\begin{lem}[Combinations of sums]
\label{lem:Combinations-of-sums} Let $S\in\N$. Then
\begin{align*}
N(S) & :=\#\big\{(a_{1},\dots,a_{k})\,\vert\,a_{j}\geq1,\sum a_{j}=S\big\}=2^{S-1}\,,\\
\tilde{N}(S) & :=\#\big\{(a_{1},\dots,a_{k})\,\vert\,a_{j}\geq1,\sum a_{j}\leq S\big\}\leq2^{S}-1\,.
\end{align*}
\end{lem}

\begin{proof}
For $S=1$, we have $N(S)=1$. Assume the claim is true for $S$.
Then for $S+1$:
\begin{align*}
\tilde{N}(S+1) & =\#\big\{(a_{1},\dots,a_{k})\,\vert\,a_{j}\geq1,\sum a_{j}\leq S+1\big\}\\
 & =\#\bigcup_{R\leq S}\big\{(a_{1},\dots,a_{k},S+1-R)\,\vert\,a_{j}\geq1,\sum a_{j}=R\big\}\cup\{(S+1)\}\\
\text{(induction)} & =1+\sum_{R\leq S}2^{R-1}=2^{S}-1
\end{align*}
On the other hand
\[
N(S+1)=\tilde{N}(S+1)-\tilde{N}(S)=2^{S+1}-2^{S}=2^{S}
\]
proves the claim.
\end{proof}
\medskip{}

We have seen in Lemma \ref{lem:Band-size-limit-by-label} that the
``size'' of a band is limited by its label $l$. To cross $n$ gaps
in our percolation model, we are more interested in the actual consecutive
stretch. It turns out that this is also just an exponential in $l$.
\begin{lem}[Total weight of a band]
\label{lem:Total-Distance-in-Band} Let $[a,b]$ be a band of label
$l$. Then,
\[
\sum_{i=a}^{b}f_{1}(i)\leq\s^{l-1}\,.
\]
\end{lem}

\begin{proof}
We have $f_{1}(i)\leq l$ for every $i\in[a,b]$. Using Lemma \ref{lem:Band-size-limit-by-label},
we have $|b-a|\leq(\tfrac{\s}{2})^{l-1}$, so
\[
\sum_{i=a}^{b}f_{1}(i)\leq l\cdot(\tfrac{\s}{2})^{l-1}\leq\s^{l-1}\,
\]
since $l\leq2^{l-1}$. 
\end{proof}
Recall from Definition \ref{def:Enumeration-Neighbours} that we may
always enumerate the ($k$) bands. The exponential decay in Lemma
\ref{lem:Exp-Decay-Band-Labels} shows that it is quite rare to encounter
bands with high labels close to the origin. This is the reason why
we may set $N_{0}:=\infty$ for a positive fraction of environments
$N=(N_{i})_{i\in\Z}$.
\begin{lem}[High labels near origin]
\label{lem:Bad-bands-near-origin}\ Consider the parameter regime
of Lemma \ref{lem:Exp-Decay-Band-Labels} for $\mathfrak{p}$ small
enough such that $24\sum_{l\geq1}(\s\mathfrak{p})^{l}<1$. Consider
the event
\[
A_{l}:=\big\{\forall\text{bands }B_{m}^{N}\text{ with }1\leq|m|\leq12\cdot\s^{l},\text{ their labels are }\leq l\big\}\,.
\]
Then, 
\[
\P(A_{l})\geq1-24\cdot(\s\mathfrak{p})^{l}\qquad\text{and}\qquad\P(\cap A_{l})>0\,,
\]
in particular, almost surely $A_{l}$ happens infinitely often.
\end{lem}

\begin{proof}
By Lemma \ref{lem:Exp-Decay-Band-Labels}, we have 
\begin{align*}
\P\big(A_{l}\big)\geq & 1-\sum_{|m|=1}^{12\cdot\s^{l}}\P(B_{m}^{N}\text{ has label }>l)\geq1-24\cdot(\s\mathfrak{p})^{l}\quad\text{and}\quad\P\big(\cap A_{l}\big)\geq1-24\sum_{l=1}^{\infty}(\s\mathfrak{p})^{l}>0\,.
\end{align*}
The last statement follows from the Borel--Cantelli lemma.
\end{proof}
\begin{proof}[Proof of Corollary \ref{cor:Adding-Infinite-Label}]
This follows from Lemma \ref{lem:Bad-bands-near-origin} and noting
that all other bands are sufficiently far away from $0$ so that they
do not merge.
\end{proof}

\subsection{Regular bands \label{subsec:Regular-Bands}}

The next point on the bucket list is making $N$ regular. $N$ being
unbounded guarantees the existence of bands of labels $\geq l$ for
all $l\in\N$ and that each such band has exactly $2$ neighbours.
We omit most proofs since they are identical to the ones in \cite{MR2116736}.
\begin{lem}[{Raising labels of maximal generators, \cite[Lemma 3.7]{MR2116736}}]
\label{lem:Raising-Maximal-Generators}\ Let $N=(N_{i})_{i\in\Z}$
be good. Let $B_{m}^{N}$ be a band of label $l$ and $i'\in\Z$ be
a maximal generator of $B_{m}^{N}$. If for all bands $B_{m'}^{N}$
of label $>l$, we have that $|m-m'|\geq\s^{l}$, then the sequence
\[
\tilde{N}_{i}=\begin{cases}
N_{i} & i\neq i'\\
N_{i}+1 & i=i'
\end{cases}
\]
satisfies the following properties:
\begin{enumerate}
\item $B_{n,k}^{N}=B_{n,k}^{\tilde{N}}\,\forall n\in\Z,k\in\N$, i.e. all
$k$ bands are identical and $\tilde{N}$ is also good.
\item If the $k$ label of $B_{n,k}^{N}$ is $t$, then the $k$ label of
$B_{n,k}^{\tilde{N}}$ is $t+\I\{i'\in B_{n,k}^{N}\}$.
\end{enumerate}
In particular, $i'$ is still a maximal generator of $B_{m}^{\tilde{N}}$.
\end{lem}

\begin{lem}[{Making $N$ more regular, \cite[Lemma 3.8]{MR2116736}}]
\label{lem:Almost-Regular-N}\ Let $N$ be good. For each $L\geq1$,
there exists $N^{L}=(N_{i}^{L})_{i\in\Z}$ such that
\begin{enumerate}
\item $N\leq N^{L}\leq N^{L+1}$,
\item $B_{m,k}^{N}=B_{m,k}^{N^{L}}$ for all $m\in\Z,k\in\N$, and
\item if $B_{m}^{N^{L}}$ and $B_{m'}^{N^{L}}$ are neighbouring bands with
label $\geq l$ and if $l\leq L$, then
\[
|m-m'|\in[\s^{l-1},\,3\cdot\s^{l-1})\,.
\]
\end{enumerate}
Furthermore, $N^{L}$ can be chosen such that $(N_{i}^{L})_{L\in\N}$
is unbounded for at most one $i$.
\end{lem}

\begin{lem}[{Making sequences regular, \cite[Lemma 3.9]{MR2116736}}]
\label{lem:Making-Sequences-Regular}\ Let $N$ be good. 
\begin{enumerate}
\item There exists a sequence $\tilde{N}\geq N$ such that all the $k$
bands for $\tilde{N}$ are identical to the $k$ bands for $N$ and
such that for neighbouring bands $B_{m},\,B_{m'}$ of label $\geq l$,
we have
\[
|m-m'|\in[\s^{l-1},\,3\cdot\s^{l-1})\,,
\]
in particular, $\tilde{N}$ is regular. In this case, we have $\tilde{N}=(N_{i})_{i\in\Z}$
with $\tilde{N}_{i}\in\N\cup\{\infty\}$ with at most one $\tilde{N}_{i}=\infty$.
(The labels may differ between $N$ and $\tilde{N}$.) 
\item There exists a sequence $\tilde{N}\geq N$ such that all the $k$
bands for $\tilde{N}$ are identical to the $k$ bands for $N$ and
such that for neighbouring bands $B_{m},\,B_{m'}$ of label $\geq l$,
we have
\[
|m-m'|\in[\s^{l-1},\,6\cdot\s^{l-1})\,,
\]
in particular, $\tilde{N}$ is regular. In this case, we have $\tilde{N}=(N_{i})_{i\in\Z}$
with $N_{i}\in\N$. (The labels may differ between $N$ and $\tilde{N}$.) 
\end{enumerate}
\end{lem}

\begin{proof}
With $N^{L}$ from Lemma \ref{lem:Almost-Regular-N}, we consider
\[
N_{i}^{\infty}:=\lim_{L\to\infty}N_{i}^{L}\in\N\cup\{\infty\}\,.
\]
We make the following observations:
\begin{enumerate}
\item If $N_{i}^{\infty}=\infty$, then $i$ must be the maximal generator
of some band $B_{m}^{N}$.
\item $N_{i}^{\infty}=\infty$ for at most one $i$. Otherwise, we would
find two separate bands $B_{m}^{N}\ni i$ and $B_{m'}^{N}\ni i'$.
The label of $B_{m}^{N^{L}}$ is bounded from below by $N_{i}^{L}$,
respectively $N_{i'}^{L}$ for $B_{m'}^{N^{L}}$. So for $l>0$ such
that $|m-m'|<\s^{l}$ and $L$ such that $\min(N_{i}^{L},\,N_{i'}^{L})\geq l$,
we would violate Lemma \ref{lem:Almost-Regular-N} Condition 3, on
the minimal distance between bands.
\end{enumerate}
Let $i^{\infty}$ be the value with $N_{i^{\infty}}^{\infty}=\infty$.
We set
\[
\tilde{N}_{i}=\begin{cases}
\lim_{L\to\infty}N_{i}^{L} & i\neq i^{\infty}\\
\infty & i=i^{\infty}
\end{cases}\,.
\]
By construction, we have that neighbouring bands $B_{m}^{\tilde{N}},\,B_{m'}^{\tilde{N}}$
always satisfy
\[
|m-m'|\in[\s^{l-1},3\cdot\s^{l-1})\,,
\]
showing the first statement. The second claim follows from choosing
$\tilde{N}_{i^{\infty}}:=N_{i}$ instead of $\infty$.
\end{proof}
\begin{rem*}[Manipulations]
 The explicit construction to make bands regular as well as Lemma
\ref{lem:Bad-bands-near-origin} allow us various manipulations on
the environment and locations of bands as well as segments.
\begin{itemize}
\item In \cite{jahnel2023continuum}, we tweak the construction such that
the origin lies not on one of the ``border segments'', but rather
on the actual inside with at least two $l$ segments distance to the
bands of label $\geq l$. Later, this ensures the existence of a circuit
around the origin. This is why we always use $12\s$ for compatibility
rather than just $6\s$.
\item In our case here, we will do quite the opposite: On a positive fraction
of environments $\Nt$, we may set $\Nt_{0}:=\infty$ without changing
any bands (Corollary \ref{cor:Adding-Infinite-Label}), effectively
considering percolation on the half-plane $\Z_{>0}\times\Z$. Ergodicity
then yields the almost-sure existence of an infinite cluster on $\Z\times\Z$. 
\end{itemize}
\end{rem*}

\subsection{Very regular bands and simple bands\label{subsec:regular-simple-bands}}

Lastly, we need a bit more information about the internal structure
of bands. This is needed to obtain crossing probabilities of strips
since we will break bands apart again. The short summary for being
very regular is: If two $k$ bands combine, then the space between
them had to be regular. The $q$ is a parameter of the distance between
those bands and will play quite an important role. 
\begin{rem*}[$k$ bands and $n$ segments]
 Short reminder that $k$ band refers to the $k$-th merging step
while $n$ segment refers to the segment between to neighbouring $(k)$
bands of label $n$.
\end{rem*}
\begin{defn}[$l$ segments (2)]
 In addition to Definition \ref{def:Enumeration-Neighbours}, we
will also call $(i_{2},i_{3})$ an $l$ \textbf{segment} if there
is a good sequence $M=(M_{i})_{i\in\Z}$ with
\[
M_{i}=N_{i}\quad\forall i\in(i_{2},i_{3})
\]
and $(i_{2},i_{3})$ is a $l$ segment for $M$. We call the segment
\textbf{regular} if it is generated by a regular sequence $M$.
\end{defn}

\begin{rem*}
The situation of the following Definition \ref{def:very-regular}
is similar to Figure \ref{fig:regular}. But since the ``neighbouring''
$n$ bands combine, they segments and bands inbetween do not have
``level'' $n-1$ but rather $q$ with $q<n-1$.
\end{rem*}
\begin{defn}[Very regular $k$ bands and $n$ segments]
\label{def:very-regular}\ Let a regular sequence $N$ be given.
\begin{enumerate}
\item Any $k$ band that is a singleton $[i,i]$ is \textbf{very regular}. 
\item The $1$ segment $\emptyset$ is \textbf{very regular}. 
\item Let $[a,d]$ be a $k$ band with label $l$ which was formed by combining
the $\tilde{k}$ bands $[a,\,b]$ and $[c,\,d]$ into the $\tilde{k}+1$
band $[a,\,d]$. $[a,\,d]$ is called \textbf{very regular} if there
are $b_{1}=b,b_{2},\dots,b_{m}$ as well as $c_{1},c_{2}\dots,c_{m-1},c_{m}=c$
with $m\leq12\s$ as well as a $q\geq1$ such that
\begin{enumerate}
\item All $\tilde{k}$ bands inside the interval $[a,\,d]$ are very regular
$\tilde{k}$ bands.
\item For all $s$, we have that $[b_{s},\,c_{s}]$ is a very regular $q$
segment.
\item For all $s<m$, we have that $[c_{s},\,b_{s+1}-1]$ is a very regular
$\tilde{k}$ band with label $q$.
\end{enumerate}
\item An $n$ segment $\mathcal{S}$ is called \textbf{very regular} if
\begin{enumerate}
\item $\mathcal{S}$ is a regular $n$ segment. (For $n=2$ and $\mathcal{S}=[a,\,b]$,
this implies $\s\leq(b-a)+2<12\s$.)
\item All $k$ bands with labels $n-1$ inside $\mathcal{S}$ are very regular.
\item All $n-1$ segments inside $\mathcal{S}$ are very regular.
\end{enumerate}
\item A band is called very regular if it is a \textbf{very regular} $k$
band for some $k$.
\item A regular sequence $N$ is called \textbf{very regular} if all the
bands generated by $N$ are very regular.
\end{enumerate}
\end{defn}

The notion of ``very regular'' allows us to split bands into smaller
parts -- enabling the induction step in Proposition \ref{prop:Drilling}.
As in Lemma \ref{lem:Making-Sequences-Regular}, we make sequences
very regular without changing the final band structure. 
\begin{lem}[{Very regular sequences, \cite[Lemma 3.12]{MR2116736}}]
\label{lem:Regular-to-Very-Regular}\ Let $N$ be good and regular.
Then, there exists $\overline{N}\geq N$ such that $\overline{N}$
is very regular and all bands and labels are identical under both
$\overline{N}$ and $N$. In particular, we may always replace a regular
sequence with a very regular sequence without changing its band structure
nor labels.
\end{lem}

\begin{proof}
This is an analogon to Lemma \ref{lem:Almost-Regular-N} and is proven
similarly (by establishing a variant of Lemma \ref{lem:Raising-Maximal-Generators}).
The labels of the final bands being unchanged follows from the construction:
To make bands very regular, one only needs to change the labels of
the $k$ bands on the ``inside''. But these labels do not contribute
to the label of the final combined band.
\end{proof}
There is one edge case that we have to worry about due to technical
issues: We want to combine bands that are close to each other first.
This led to the quite cumbersome merging scheme in Definition \ref{def:bands-and-labels}/Algorithm
\ref{alg:Merging-Indices} as well as the following:

\begin{defn}[Simple $k$ bands]
\label{def:simple-bands} \ 
\begin{enumerate}
\item Any $k$ band that is a singleton $[i,i]$ is \textbf{simple}. 
\item Let $[a,d]$ be a $k$ band with label $l$ which was formed by combining
the $\tilde{k}$ bands $[a,\,b]$ and $[c,\,d]$ into the $\tilde{k}+1$
band $[a,\,d]$. $[a,\,d]$ is called \textbf{simple} if both $[a,b]$
and $[c,d]$ are simple $\tilde{k}$ bands as well as 
\[
1+D_{\tilde{k}}(b,c)<(12\s)^{2}
\]
(see Definition \ref{def:bands-and-labels}, Algorithm \ref{alg:Merging-Indices}).
\end{enumerate}
\end{defn}

\begin{rem*}[$q$ in simple bands]
 By Definition \ref{def:simple-bands} and Algorithm \ref{alg:Merging-Indices},
we see that simple bands satisfy $q\leq2$ with $q$ as in Definition
\ref{def:very-regular} above. Furthermore, if 
\[
\s\geq72=(12)^{2}/2\,,
\]
then this is even an equivalence since for $q=3$, we would automatically
have
\[
1+D_{\tilde{k}}(b,c)\geq2\cdot\s^{3}>(12\s)^{2}\,.
\]
(Using that the minimal size of a $3$ segment is $\s^{2}$.) This
allows for an easy characterisation.
\end{rem*}
\begin{lem}[Sufficient criterion for simple bands]
\label{lem:sufficient-criterion-simple-band} Let $[a,d]$ be a $k$
band with label $l$ which was formed by combining the $\tilde{k}$
bands $[a,\,b]$ and $[c,\,d]$ into the $\tilde{k}+1$ band $[a,\,d]$.
If 
\[
1+D_{\tilde{k}}(b,c)<(12\s)^{2}\,,
\]
then $[a,b]$ and $[c,d]$ also had to be simple $\tilde{k}$ bands.
In particular, if $\s\geq72$ and $q\leq2$, then $[a,d]$ is simple.
\end{lem}

\begin{proof}
One checks that if either $[a,b]$ or $[c,d]$ have been non-simple,
then it would contradict with Step 2 in the construction of $k$ bands
in Definition \ref{def:bands-and-labels}. The last statement follows
from the previous remark.
\end{proof}
The nice thing about simple bands -- and the sole reason we need
to look at them -- is that their ``stretch'' grows at most linearly
in $l$ (rather than the extremely crude exponential estimate in Lemma
\ref{lem:Total-Distance-in-Band}):
\begin{lem}[Maximal stretch of simple bands]
\label{lem:maximal-stretch-of-simple} Let $\s\geq72$. Let $[a,d]$
be a simple $k$ band with label $l\geq2$. Then,
\[
\sum_{i\in[a,d]}f_{k}(i)\leq l+(13\s)^{2}\cdot(l-2)/2\,.
\]
\end{lem}

\begin{proof}
In the case of a singleton $[a,d]=[a,a]$, this is true since $f_{k}(a)=l$.
Now, assume that the claim is true for all $k$ bands with labels
$<l$. If the $k$ band is not a singleton, we split it up into the
simple $\tilde{k}$ bands $[a,\,b]$ and $[c,\,d]$ as before with
labels $l_{1}$ respectively $l_{2}$, where $l_{1}+l_{2}=l$. Since
$[a,d]$ is simple and $\s\geq72$, it is very regular with $q\leq2$.
Therefore, there can at most be $12\s$ bands of label $2$ between
$[a,b]$ and $[c,d]$ with the rest being bands of label $1$. Now
by the induction hypothesis
\begin{align*}
\sum_{i\in[a,d]}f_{k}(i) & =\sum_{i\in[a,b]}f_{\tilde{k}}(i)+\sum_{i\in[c,d]}f_{\tilde{k}}(i)+\sum_{i\in(b,c)}f_{k}(i)\\
 & \leq\big\{ l_{1}+(13\s)^{2}\cdot(l_{1}-2)/2\big\}+\big\{ l_{2}+(13\s)^{2}\cdot(l_{2}-2)/2\big\}+\big\{2\cdot12\s+(12\s)^{2}\big\}\\
 & \leq l+(13\s)^{2}\cdot(l-4)/2+(13\s)^{2}=l+(13\s)^{2}\cdot(l-2)/2\,,
\end{align*}
which shows the claim. 
\end{proof}

We conclude the section with parameter estimates on very regular bands.
These turn out to be quite crucial, in particular the upper bound
for $q$.

\begin{lem}[Estimates for $m,r,q$ on very regular bands]
\label{lem:q-estimates} Assume that we have split the very regular
band into bands with labels $m,r$ and have the space inbetween with
parameter $q$. Then,
\begin{equation}
m+r=n\quad\forall q\leq8\label{eq:m-plus-r-equals-n}
\end{equation}
as well as
\begin{equation}
m+r-\lfloor\cD q\rfloor=n+\sigma\label{eq:m-plus-r-minus-q-is-n}
\end{equation}
with $\sigma\in\{-1,0,1\}.$ Furthermore, 
\begin{equation}
q\leq\lfloor(2-\cD)^{-1}n\rfloor=:\flat(n)\,.\label{eq:q-smaller-flat-n}
\end{equation}
Note that since $\cD<1/11$, we have $\flat(n)\leq\lfloor\tfrac{11}{21}n\rfloor$.
\end{lem}

\begin{proof}
We get to return to the label generation again (Definition \ref{def:bands-and-labels}):
\[
n=m+r-\lfloor\cD\log_{\s}(1+D)\rfloor
\]
where $D$ is the number of bands between the bands of label $m,r$
right before combining. Since the bands are very regular, we have
at most $12\cdot\s-1$ many bands of label $q$ between them with
corresponding $q$ segments. Each $q$ segment contains at least $\s^{q-1}$
and at most $12\cdot\s^{q-1}$ many bands. Therefore
\begin{align*}
\s^{q-1} & \leq D\leq12\cdot\s^{q-1}\cdot12\s+(12\s-1)\\
\s^{q-1} & \leq1+D\leq\s^{q}\cdot13^{2}\\
q-1 & \leq\log_{\s}\big(1+D\big)\leq q+2\log_{\s}13\\
\cD q-\cD & \leq\cD\log_{\s}\big(1+D\big)\leq\cD q+\cD2\log_{\s}13\,.
\end{align*}
If $q\leq8$, then $\cD q+\cD2\log_{\s}13<1$, in particular $\lfloor\cD\log_{\s}(1+D)\rfloor=0$.
This proves Equation (\ref{eq:m-plus-r-equals-n}). Furthermore, since
$\cD(2\log_{\s}13+1)\leq\tfrac{1}{11}(1+1)<1$, we have
\[
\big|\lfloor\cD q\rfloor-\lfloor\cD\log_{\s}\big(1+D\big)\rfloor\big|\leq1
\]
which yields Equation (\ref{eq:m-plus-r-minus-q-is-n}). Since $m,r>q$,
we have
\begin{align*}
n+\text{either }0\text{ or }1\geq & 2q-\lfloor\cD q\rfloor+2\\
n\geq & 2q-\cD q=q(2-\cD)\\
(2-\cD)^{-1}n\geq & q\,,
\end{align*}
i.e., Inequality (\ref{eq:q-smaller-flat-n}) since $q$ is an integer.
\end{proof}
\begin{rem*}[Final remarks]
 As alluded to early on, we will use the whole ``segment-band''
framework for both the temporal rows as well as spatial columns. In
the case of the spatial columns, we will attempt to cross bad bands
in a single jump, so not much of the inner structure is needed.

The temporal columns are much harder to handle. We will need to exploit
that bands are very regular in order use induction. Lemma \ref{lem:q-estimates}
will also play a crucial role throughout Section \ref{subsec:Drilling}
as it limits us in how thin we can make strips. The notion of simple
bands is needed for the base case of $q\leq2$.
\end{rem*}

\section{Details: proving percolation \label{sec:Proving-Percolation}}

We employ the band/segment grouping scheme for the time/space stretches
$(\Nt_{t})_{t\in\Z}$, $(\Nx_{x})_{x\in\Z}$ with parameters $\s_{t},\s_{x}$
and $\cD=1/12$. We may assume without loss of generality that these
stretches are very regular (Lemma \ref{lem:Making-Sequences-Regular},
\ref{lem:Regular-to-Very-Regular}).

\subsection{Connectivity inside/between good boxes}

The usual idea with multiscale/block arguments is to connect boxes
of different levels with each other. Directionality adds bloat to
the proofs, but the principle behind is actually simple and graphical: 

\begin{figure}
\includegraphics[width=1\columnwidth]{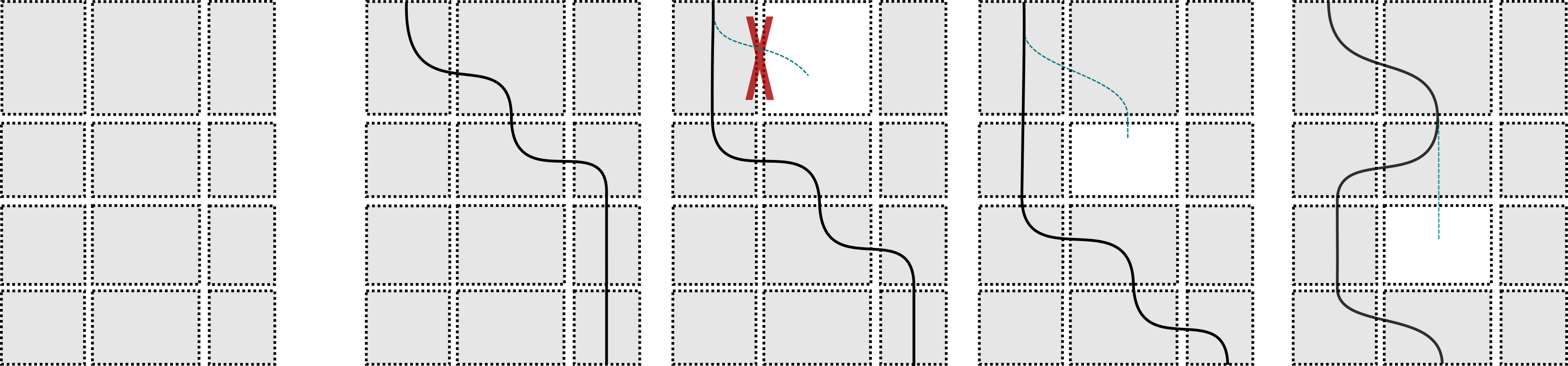}\caption{Depicted are the schemes by which we connect boxes. In the best case,
we just go to the target column and move straight down. Otherwise,
we have to dodge bad boxes/connections. \label{fig:connecting-boxes}}
\end{figure}

\begin{lem}[Reachable boxes]
\label{lem:Reachable-Boxes} Let a rectangular area of $l_{x}\geq2$
columns and $l_{t}$ rows of $n$ boxes be given, which are separated
by $n$ gaps and $(n+1,n)$ strips. Number them by 
\[
(B_{i,j})_{1\leq i\leq l_{t},1\leq j\leq l_{x}}\,.
\]
Assume that Lemma \ref{lem:Connecting-Inputs-Outputs-Inside} is true
for $n$. If at most one of the $n$ boxes, $n$ gaps or $(n+1,n)$
strips is bad, then for any good $n$ boxes $B_{i,j}$ and $B_{i',j'}$
with $(i'-i)\geq l_{x}$, we have
\[
\Inv(B_{i,j})\fc\Outv(B_{i',j'})\,.
\]
\end{lem}

\begin{proof}
Without loss of generality, we assume that $j\leq j'$, otherwise
we mirror the whole procedure. We sketch the connecting procedure
in Figure \ref{fig:connecting-boxes} where horizontal connections
are made via 
\[
\Inv(B_{k,l})\fc\Outh(B_{k,l})\leadsto\Inh(B_{k,l+1})\fc\Outv(B_{k,l+1})
\]
and vertical ones via $\Out(B_{k,l})\leadsto\In(B_{k+1,l})$. First
of all, it suffices to only look at the case with at most one bad
box: If the gap between $B_{k,l}$ and $B_{k,l+1}$ is bad, we simply
declare $B_{k,l}$ to be bad (if it is not the starter box, otherwise
take $B_{k,l+1}$). The same works for strips. We distinguish two
cases.
\begin{enumerate}
\item The procedure is straight-forward. If we are currently in $\Inv(B_{k,l})$
with $l<j'$ and \textbf{both} $B_{k,l+1}$, $B_{k+1,l+1}$ are good,
then move towards $\Outv(B_{k,l+1})$ (and then $\Inv(B_{k+1,l+1})$).
Otherwise, simply move down $\Inv(B_{k,l})\leadsto\Inv(B_{k+1,l})$
and proceed.
\item If $l=j'$ -- i.e.\ we have reached the target column $j'$ --
and $B_{k+1,l}$ is good, then again $\Inv(B_{k,l})\leadsto\Inv(B_{k+1,l})$.
Otherwise, dodge to side, i.e.\ $\Inv(B_{k,l})\leadsto\Inv(B_{k+1,l+\sigma})$
with $\sigma=1$ if $l=1$ and $-1$ otherwise.
\end{enumerate}
Since at most one box is bad, only one ``delaying case'' can happen,
so we still reach the target output.
\end{proof}
Now, we know that rectangular regions of good $n$ boxes are well-connected,
given that Lemma \ref{lem:Connecting-Inputs-Outputs-Inside} holds.
Naturally, we have to prove said lemma now for $n=1$ and $n+1$.
\begin{proof}[Proof of Lemma \ref{lem:Connecting-Inputs-Outputs-Inside}]
 This is true for $n=1$. For $n+1$, the claim on $B_{n}$ follows
from Lemma \ref{lem:Reachable-Boxes} and Equation (\ref{eq:kappa_ver-and-kappa_hor})
on $\s_{x},\s_{t}$: We again number the boxes in $B_{n}$ as $(B_{i,j})_{1\leq i\leq l_{t},1\leq j\leq l_{x}}$
and take some $v\in\Inv(B_{n})$, in particular $v\in\Inv(B_{1,j})$.
If $w\in\Outh(B_{n})$, then in particular $w\in\Outh(B_{i,j'})$
for some $B_{i,j'}$. Furthermore, both $B_{i,j'}$ and $B_{i-1,j'}$
are good (by definition of $\Outh(B_{n})$) with 
\[
i-1\geq12\s_{x}+1\geq l_{x}\,.
\]
Using Lemma \ref{lem:Reachable-Boxes} and the induction hypothesis,
we connect $v\fc\Outv(B_{i-1,j'})\leadsto\Inv(B_{i,j'})\fc\Outh(B_{i,j'})$.
This proves the first part. The second part follows directly from
Lemma \ref{lem:Reachable-Boxes} and the third part works analogously
to the first part. Finally,
\[
A\fc B\leadsto C\fc D\implies A\fc D\,,
\]
yields the final statement on neighbouring boxes.
\end{proof}

\subsection{Connecting outputs with inputs and multiscale estimates}

Not all outputs of $n$ boxes connect directly to inputs. There is
always some loss due to bad $n-1$ boxes in prior steps. In this subsection,
we quantify the minimum amount of suitable connectors, which yields
the probability of good $(n+1,n)$ strips as well as $n$ gaps.
\begin{defn}[$(\kver,n)$ trees]
~\label{def:kappa-trees}
\begin{enumerate}
\item A $(\kver,1)$ \textbf{tree} is any single vertex.
\item A $(\kver,n)$ \textbf{tree} consists of $\kver$ many disjoint $(\kver,n-1)$
trees such that they all lie inside $\{t\}\times(x_{1},x_{2}]$ for
some some $t\in\Z$ and spatial $n$ segment $(x_{1},x_{2})$.
\end{enumerate}
\end{defn}

\begin{rem*}
$(\kver,n)$ trees capture the basic shape of the sets $\Inv(B_{n})$
and $\Outv(B_{n})$ of an $n$ box $B_{n}$. Each such tree contains
(exactly) $\kver^{n-1}$ many vertices. They will play the role of
``connectors'' between vertically neighbouring boxes as we see in
the following:
\end{rem*}
\begin{lem}[$(\kver,n)$ trees between good $n$ boxes]
\label{lem:Trees-between-boxes} Let $B_{n}$ and $B_{n}'$ be $n$
boxes where $B_{n}$ lies on top of $B_{n}'$ (only separated by an
$(n+1,n)$ strip). Then, they define at least one $(\kver,n)$ tree
$T$ such that:
\[
T\subset\Outv(B_{n})\qquad\text{and}\qquad\pi_{x}(T)\subset\pi_{x}(\Inv(B_{n}'))
\]
where $\pi_{x}:\Z^{2}\to\Z$ is the projection onto the $x$-coordinate.
In words: They define a $(\kver,n)$ tree $T$ such that $T$ lies
in the same column as $\Out(B_{n})$ and $\In(B_{n}')$. 
\end{lem}

\begin{proof}
The proof is by induction. In the case a $1$ box $B_{n}=[t_{1},t_{2}]\times\{x\}$
we take $T=\{(t_{2},x)\}$. For general $n$, we know that in each
row of these at least $\kver+2$ many $n-1$ boxes, at most one of
these boxes is bad. Therefore, there are at least $\kver$ many pairs
of good vertically neighboured $n-1$ boxes. By the induction hypothesis,
these define $\kver$ many $(\kver,m-1)$ trees which satisfy Condition
2 of Definition \ref{def:kappa-trees} above since they lie in $B_{n}$.
Therefore, we obtain a $(\kver,n)$ tree as claimed.
\end{proof}
This covers the case of vertical connectors. We set up the same framework
analogously for horizontal connectors, but actually keep things straight
and explicit here:
\begin{lem}[Number of horizontal connectors between good $n$ boxes]
\label{lem:Number-Horizontal-Connectors} Let $B_{n},B_{n}'$ be
neighbouring good $n$ boxes. Then, there are at least $\khor^{n}$
many edges from $\Outh(B_{n})$ to $\Inh(B_{n}')$ crossing exactly
over the $n$ gap inbetween..
\end{lem}

\begin{proof}
In the case of $1$ boxes $B_{1}=[t_{1},t_{2}]\times\{x\},\,B_{1}'=[t_{1},t_{2}]\times\{x'\}$,
every $(t,x)\in\Outh(B_{1})$ has an outgoing edge to $(t+1,x')\in\Inh(B_{2})$
for every $t\in[t_{1},t_{2})$. This makes $|t_{2}-t_{1}|\geq\lceil\s_{t}/12\rceil-1\geq\khor$
many different edges. 

For the case of the $n+1$ boxes $B_{n+1},B_{n+1}'$, we see by the
definition of inputs/outputs (in Definition \ref{def:Good-boxes})
that $\Outh(B_{n+1})$ and $\Inh(B_{n+1}')$ consist of $\khor+4$
many opposing $n$ boxes if they were all valid. Since $B_{n+1},B_{n+1}'$
are good, at most $2$ of the the boxes in $\Outh(B_{n+1})$ might
not be valid, same for $\Inh(B_{n+1}')$. Therefore, we have $\khor$
many opposing $n$ boxes that may connect with each other. By the
induction hypothesis, each of these contribute at least $\khor^{n}$
many edges, so we have $\khor\cdot\khor^{n}=\khor^{n+1}$ in total
which proves the claim.
\end{proof}
With this, we have guaranteed that there are exponentially many potential
connectors for both the vertical strips as well as horizontal gaps.
This is important since we want to use the following estimate: 
\begin{lem}[Combinatorial estimate]
\label{lem:Multiscale-Estimate} Assume there is a collection of
at most $C$ ``objects'' that are each good with probability at
least $P_{n}$ independently from each other. Furthermore, assume
that a certain object of level $n+1$ is good if at most one of the
$C$ prior objects is bad. Then, for any $\pp\in(0,1)$ with $\pp^{n+1}\leq C^{-6}$,
if $n\geq1$ and
\[
1-P_{n}\leq\pp^{n+1}\,,
\]
then also
\[
1-P_{n+1}\leq\pp^{n+2}\,.
\]
\end{lem}

\begin{proof}
We first write $1+k_{n}:=(1-P_{n})^{-1}$. The level $n+1$ object
is good with probability at least
\[
P_{n+1}\geq(P_{n})^{C}+C\cdot(P_{n})^{C-1}\cdot(1-P_{n})\,.
\]
Therefore
\[
1-P_{n+1}\leq1-\left[\left(\frac{k_{n}}{1+k_{n}}\right)^{C}+C\cdot\left(\frac{k_{n}}{1+k_{n}}\right)^{C-1}\cdot\frac{1}{1+k_{n}}\right]=\frac{(k_{n}+1)^{C}-(k_{n})^{C}-C\cdot(k_{n})^{C-1}}{(1+k_{n})^{C}}\,.
\]
The subtrahends are exactly the first two terms in this binomial expression.
Therefore, 
\begin{align*}
1-P_{n+1}=(1+ & k_{n})^{-C}\cdot\sum_{i=0}^{C-2}{C \choose i}(k_{n})^{i}\leq\frac{C^{3}\cdot(1+k_{n})^{C-2}}{(1+k_{n})^{C}}\\
 & \leq(1+k_{n})^{-1.5}\leq(1-P_{n})^{1.5}\leq\pp^{(n+1)\cdot1.5}\leq\pp^{n+2}\,,
\end{align*}
where we also used $1+k_{n}=(1-P_{n})^{-1}\geq\pp^{-(n+1)}\geq C^{6}$.
\end{proof}
\begin{rem*}
The lemma can be generalised to allow for $C[n]=a_{1}\e^{a_{2}n}$
instead of a constant $C$.
\end{rem*}
In our case, the ``level $n+1$'' object will be an $n+1$ box containing
up to $C$ many $n$ boxes, $(n+1,n)$ strips as well as $n$ gaps
inbetween. By construction, each $n+1$ box will then contain at most
$(12\s_{x}+1)\cdot12\s_{t}$ many $n$ boxes, so the total number
of level $n$ objects is
\begin{equation}
C\leq(12\s_{x}+1)\cdot12\s_{t}\cdot(1+1+1)\leq450\cdot\s_{x}\s_{t}\,.\label{eq:Number-Objects-inside-Box}
\end{equation}
There is a small technical issue in using Lemma \ref{lem:Multiscale-Estimate}:
In order to ensure a high probability for $n$ gap crossings, we need
a large amount of connectors, i.e.\ $\s_{t}$ to be large. But this
also results in a larger constant $C$, so the gap crossing probability
has to grow accordingly. The next two lemmas ensure that this circular
dependency is not a problem.
\begin{lem}[Horizontal strip crossing]
\ Let $B_{n}$ and $B_{n}^{'}$ be neighbouring good $n$ boxes.
Then 
\[
\P\big\{\Outh(B_{n})\not\leadsto\Inh(B_{n}')\big\}\leq\exp\big(-\{(1+\s_{x})^{-\alpha}\khor\}^{n}\big)\,.
\]
\end{lem}

\begin{proof}
By Lemma \ref{lem:Number-Horizontal-Connectors}, there are at least
$\khor^{n}$ suitable edges that would connect $\Outh(B_{n})$ with
$\Inh(B_{n}')$ if they were open. By Lemma \ref{lem:Total-Distance-in-Band},
these edges have length at most $\s_{x}^{n}$. Therefore
\begin{align*}
\P\big\{\Outh(B_{n}) & \not\leadsto\Inh(B_{n}')\big\}\leq\big(1-\{1+\s_{x}^{n}\}^{-\alpha}\big)^{(\khor^{n})}\leq\exp\big(-\{1+\s_{x}\}^{-n\alpha}\cdot\khor^{n}\big)\\
\leq & \exp\big(-\{(1+\s_{x})^{-\alpha}\khor\}^{n}\big)\,,
\end{align*}
which yields the claim.
\end{proof}
\begin{lem}[Ensuring high probability of horizontal strip crossings]
\label{lem:Gap-Crossing}\ Given fixed $\pp$, $\s_{x}$ and $\alpha$,
then we have for $\s_{t}$ large enough (equivalently $\khor$ large
enough): For any $n$ gap $G$, we have
\[
1-\P(G\text{ is good})\leq\min\left\{ \pp^{n+1},\left(450\cdot\s_{x}\s_{t}\right)^{-6}\right\} \,.
\]
In particular, we may ensure that both Theorem \ref{lem:Main-Thm}
Point 2 as well as the requirements of Lemma \ref{lem:Multiscale-Estimate}
hold for horizontal gaps.
\end{lem}

\begin{proof}
Using the previous lemma, we see that we only need to show
\[
2\exp\big(-\{(1+\s_{x})^{-\alpha}\khor\}^{n}\big)\leq\min\left\{ \pp^{n+1},\left(450\cdot\s_{x}\s_{t}\right)^{-6}\right\} \,.
\]
First, by Equation (\ref{eq:kappa_ver-and-kappa_hor})
\[
\khor\geq\s_{t}/12-25\s_{x}=\s_{t}/12-c\,,
\]
so the requirements on horizontal crossings are met if both
\[
2\exp\big(-\big\{(1+\s_{x})^{-\alpha}(\s_{t}/12-c)\big\}^{n}\big)\leq\exp\big(-\big\{(1+\s_{x})^{-\alpha}\khor\big\}^{n}\big)\leq\left(450\cdot\s_{x}\s_{t}\right)^{-6}
\]
and
\begin{align*}
\exp\big(-\big\{(1+\s_{x})^{-\alpha}\khor\big\}^{n}\big) & \leq\pp^{n+1}\\
\big\{(1+\s_{x})^{-\alpha}\khor\big\}^{n} & \geq(n+1)\log\tfrac{1}{\pp}
\end{align*}
are satisfied, which is true for $\s_{t}$ (equivalently $\khor)$
large enough.
\end{proof}

\subsection{Proof of Lemma \ref{lem:Main-Thm}}

\begin{figure}
\includegraphics[width=1\columnwidth]{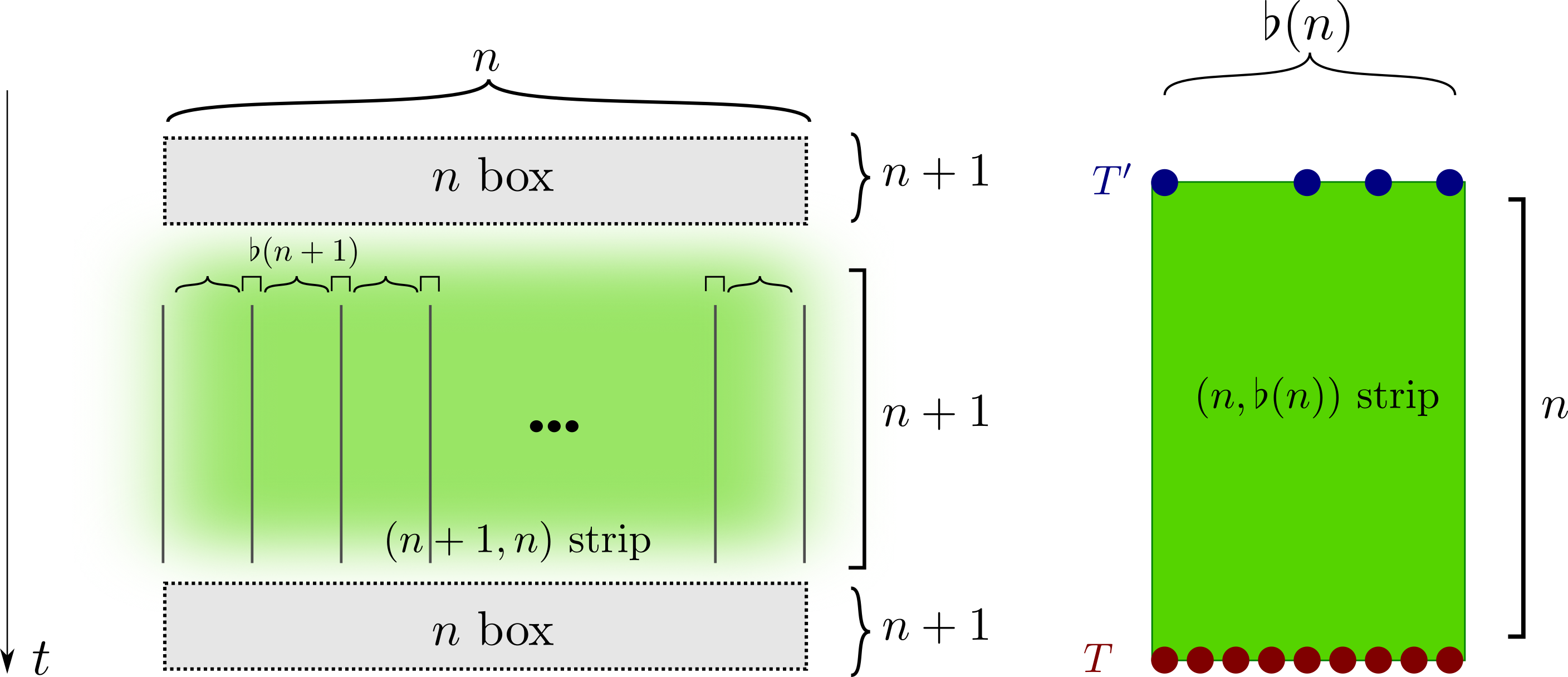}\caption{As always, curly brackets indicate bands with the square brackets
indicating bands. Whenever we consider a crossing of a $(n+1,n)$
strip, we actually try to do so in disjoint $(n+1,\flat(n+1))$ strips.
Since these are \textquotedblleft thin\textquotedblright{} objects,
they may be broken down further, so we pay special attention to $(n,\flat(n))$
strips.\label{fig:main-thm-setting}}
\end{figure}

We can now prove Lemma \ref{lem:Main-Thm} provided that Proposition
\ref{prop:Drilling} below holds for $n+1$. The setting is depicted
in Figure \ref{fig:main-thm-setting}.
\begin{prop}[Drilling]
\label{prop:Drilling} Let $S$ be a $(n,\flat(n))$ strip with $n\geq1$.
Let $T'$ be a collection of $(\kver,k)$ trees on top of $S$ and
$T$ be a $(\kver,\flat(n))$ tree on the bottom of $S$ with $\pi_{x}(T')\subset\pi_{x}(T)$
where $\pi_{x}$ is the projection onto the $x$-coordinate. Then,
\[
\P(\exists\text{ a crossing of }S\text{ intersecting both }T\text{ and }T')\geq\kver^{-\flat(n)}\cdot\#T'\,.
\]
\end{prop}

\begin{proof}[Proof of Theorem \ref{lem:Main-Thm}]
 Using Lemma \ref{lem:Gap-Crossing}, Point $2$ is ensured by fixing
some large $\khor$ (or equivalently $\s_{t})$. WLOG, we assume $\pp\leq(450\cdot\s_{x}\s_{t})^{-6}$
(for Lemma \ref{lem:Multiscale-Estimate}). Then, we choose $p$ large
enough such that Lemma \ref{lem:Main-Thm} holds for every $n\leq\mathcal{N}$,
where $\mathcal{N}$ comes from Equation (\ref{eq:Induction-Start-N})
below. We also require $p^{100\s_{t}{}^{2}}(1-e^{-1})\geq\kver^{-1/2}$
in Equation (\ref{eq:induction-p}). 
\begin{enumerate}
\item[3)]  We show that Point 3 holds for $n$ given that Proposition \ref{prop:Drilling}
holds for $n+1$. Let $\mathcal{N}\in\N$ large enough such that 
\begin{equation}
\kver^{n-\flat(n)-2}\geq(n+1)\log\tfrac{1}{\pp}\label{eq:Induction-Start-N}
\end{equation}
for every $n\geq\mathcal{N}$. We use Lemma \ref{lem:Trees-between-boxes}
to first get a $(\kver,n)$ tree $\tilde{T}\subset\Outv(B_{n})$ with
$\pi_{x}(\tilde{T})\subset\pi_{x}(\Inv(B_{n}'))$. Now, the $(n+1,n)$
strip can be divided into $(2+\kver)^{n-\flat(n+1)}$ many $(n+1,\flat(n+1))$
strips. We will choose (exactly) $\kver^{n-\{\flat(n)+1\}}$ disjoint
$(n+1,\flat(n+1))$ strips $S$ such that they have a $(\kver,\flat(n+1))$
tree $T'$ on top satisfying $T'\subset\tilde{T}$. By Proposition
\ref{prop:Drilling}, the probability of crossing $S$ is at least
$\kver^{-\flat(n+1)}\cdot\#T'=1/\kver$. Since all those strips are
disjoint, these events are independent. Therefore, we have
\begin{align*}
\P\big\{\nexists\text{ a crossing of } & \bar{S}\text{ intersecting }\Outv(B_{n}),\Inv(B_{n}')\big\}\\
\leq & \left(1-1/\kver\right)^{\kver^{n-\flat(n+1)}}\leq\exp\big(-\kver^{n-\flat(n)-2}\big)\leq\pp{}^{n+1},
\end{align*}
which shows Point 3.
\item[1)]  Showing Point 1 for $n+1$ is a straight-forward application of
Lemma \ref{lem:Multiscale-Estimate} after using all the estimates
on $n$ boxes, $(n+1,n)$ strips and $n$ gaps.
\end{enumerate}
\end{proof}
Judging by the remaining pages, one can guess that Proposition \ref{prop:Drilling},
i.e.,\textbf{ drilling}, is the most difficult part. Also the fact
that we have yet to use that $\Nt$ is very regular. The good news
is that we can already prove the case of simple bands.
\begin{proof}[Proof of Proposition \ref{prop:Drilling} for simple bands]
 The case of simple bands is equivalent to $q\leq2$ (see Lemma \ref{lem:sufficient-criterion-simple-band}).
We assume that the temporal $n$ band generating the $(n,k)$ strip
is simple with $k\geq n/2$. We generate crossings by going straight
through a column. By Lemma \ref{lem:maximal-stretch-of-simple} (and
using that $\s_{t}>17'000$), this probability is at least
\[
p^{n+(13\s_{t})^{2}(n-2)/2}\geq p^{100\s_{t}{}^{2}n}.
\]
There are $\#T'$ vertices (or rather columns) which potentially form
an appropriate crossing if they were open. Thus, using our assumption
of 
\begin{equation}
p^{100\s_{t}{}^{2}n}(1-e^{-1})\geq\khor^{-n/2}\geq\khor^{-k}\label{eq:induction-p}
\end{equation}
 as well as Lemma \ref{lem:Technical-Estimate} below
\begin{align*}
\P(\exists\text{ a cluster in }S & \text{ connecting }T\text{ and }T')\geq1-(1-p^{100\s_{t}{}^{2}n})^{\#T'}\\
\geq & \min\left\{ 1-e^{-1},\,\#T'\cdot p^{100\s_{t}{}^{2}n}(1-e^{-1})\right\} \geq\kver^{-k}\#T'\,,
\end{align*}
which proves the case of simple bands. (Note that $\#T'\leq\kver^{k-1}$.)
\end{proof}
Here is the auxiliary lemma we previously used and will continue to
use in the future.
\begin{lem}[{\cite[Lemma 4.2]{MR2116736}}]
\label{lem:Technical-Estimate}\ For any $c,p_{1},\dots,p_{n}$
with $0<p_{i}<1$ and $a:=\sum_{1}^{n}p_{i}$, we have
\[
1-\prod_{i=1}^{n}(1-p_{i})\geq\min\left\{ 1-e^{-c},\,\tfrac{a}{c}(1-e^{-c})\right\} \,.
\]
\end{lem}

\subsection{Drilling: preparation}

\begin{figure}

\includegraphics[width=1\columnwidth]{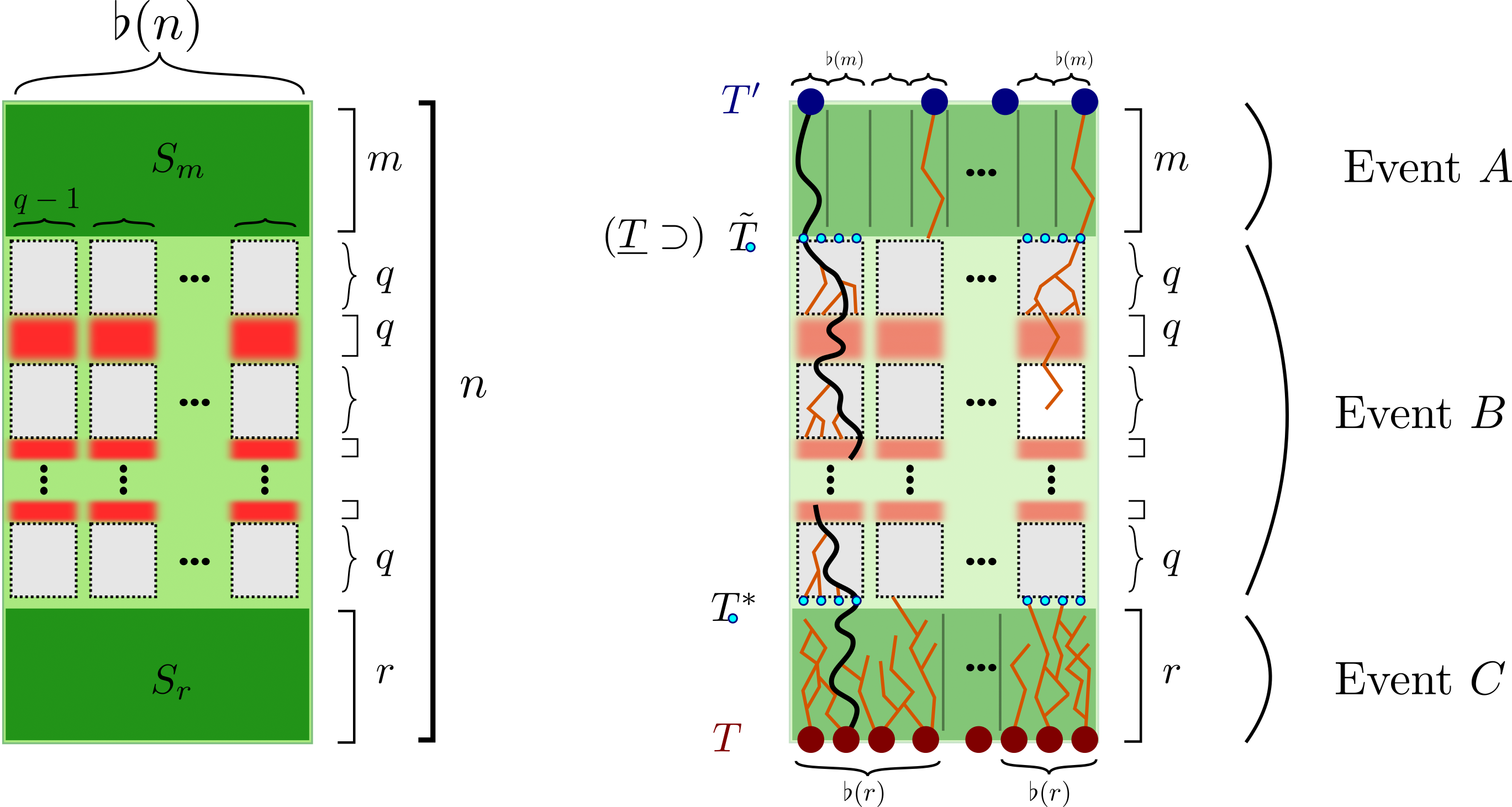}\caption{Drilling/generating a crossing of $S$. The Events $A,B,C$ together
yield the crossing (bold black path) with a minimal probability depending
on $\#T'$.\label{fig:strip-structure}}

\end{figure}

Now comes the tough part. Assume that Lemma \ref{lem:Main-Thm} holds
until $\flat(n)\leq n-2$. We want to see that we can \textbf{drill}
through arbitrary $(n,\flat(n))$ strips $S$, i.e., Proposition \ref{prop:Drilling}
holding even for $q\geq3$. We will use that the temporal stretches
$\Nt$ are very regular to break up $S$ into three smaller parts,
see Figure \ref{fig:strip-structure} with the other variables being
introduced during the course of this section. On the top, we have
a $(m,\flat(n))$ strip $S_{m}$. On the bottom, we have a $(r,\flat(n))$
strip $S_{r}$. In the middle, there are up to $12\s_{t}$ rows of
$q-1$ boxes separated by $(q,q-1)$ strips. Lemma \ref{lem:q-estimates}
will be crucial in our endeavour.

The outline of the remaining proof is as follows. If
\begin{enumerate}
\item[(A)] there are ``enough'' crossings of $S_{m}$ which intersect $T'$
(Equation (\ref{eq:Enough}), Lemma \ref{lem:sufficiently-many-crossings}),
\item[(B)] these crossings survive through the column of $q-1$ boxes to $S_{r}$
(Lemma \ref{lem:Good-Columns}),
\item[(C)] one of these survivors connects in $S_{r}$ to $T$ (Proposition
\ref{prop:Drilling}),
\end{enumerate}
then there exists a crossing of $S$ intersecting $T'$ and $T$.
For Event $B$, a single crossing survives with probability at least
$0.99$ (Lemma \ref{lem:Lemma-46}). This is a rather simple calculation.
As for the rest, the technicalities are more difficult than the actual
proof.

In Lemma \ref{lem:pooling-strips}, we pool together small strips
and estimate the probability of a crossing happening for at least
one of them. Then, we estimate the probability of Event $A$ in Lemma
\ref{lem:sufficiently-many-crossings}. We do so by pooling several
$(m,\flat(m))$ strips together so that each such collection has a
sufficiently high probability of crossing $S_{m}$. Lemma \ref{lem:pooling-strips}
allows us to pool together all survivors from Event $A$ to obtain
a lower bound on the probability of Event $C$. Finally, Proposition
\ref{prop:Drilling} follows from combining all of the previous calculations.

\medskip{}

Let us briefly consider a general $(j,J)$ strip $S^{*}$ with $j<n$
and $J\geq\flat(j)$. Let $\hat{S}:=\cup\hat{S}_{i}$ be a disjoint
union of $(j,\flat(j))$ strips in $S^{*}$. Let $T^{*}$ (the target)
be a $(\kver,J)$-tree on the bottom of $S^{*}$ which intersects
each $\hat{S}_{i}$ in a $(\kver,\flat(j))$-tree. Let $\hat{T}$
be a union of $l$-trees on top of $\hat{S}$ where $l\leq\flat(j)$,
all lying in the columns of $T^{*}$.
\begin{lem}[{Pooling together strips for crossings, \cite[Lemma 4.4]{MR2116736}}]
\label{lem:pooling-strips}\ Suppose Proposition \ref{prop:Drilling}
holds for $j\leq n-2$. Then,
\begin{align*}
\P(\exists\text{ a crossing of }\hat{S}\text{ intersecting }\hat{T}\text{ and }T^{*} & )\geq\min\left\{ 0.9,\,\tfrac{1}{3}\kver^{\flat(j)}\cdot\#\hat{T}\right\} \,.
\end{align*}
Each such crossing is confined to its respective $(j,\flat(j))$ strip.
\end{lem}

\begin{proof}
$\hat{T}$ is a union of $(\kver,l)$ trees. Let $\hat{T}=\cup\hat{T}_{i}$
where $\hat{T}_{i}$ consists of the $(\kver,l)$ trees belonging
to $\hat{T}$ that lie inside the $(j,\flat(j))$ strip $\hat{S}_{i}$
(recall $l\leq\flat(j)$). By the induction hypothesis, we have 
\[
\P(\exists\text{ a crossing of }\hat{S_{i}}\text{ intersecting }\hat{T_{i}}\text{ and }T^{*})\geq\kver^{-\flat(j)}\cdot\#\hat{T_{i}}\,.
\]
These are independent events since all the $\hat{S}_{i}$ are disjoint.
Lemma \ref{lem:Technical-Estimate} with $c=2.31$ yields
\begin{align*}
\P(\exists\text{ a crossing of } & \hat{S}\text{ intersecting }\hat{T}\text{ and }T^{*})\\
 & \geq1-\prod_{i}\big(1-\P\{\exists\text{ a crossing of }\hat{S_{i}}\text{ intersecting }\hat{T_{i}}\text{ and }T^{*}\}\big)\\
 & \geq(1-e^{-2.31})\min\Big\{1,\,\tfrac{1}{2.31}\kver^{-\flat(j)}\sum_{i}\#\hat{T}_{i}\Big\}\geq\min\left\{ 0.9,\,\tfrac{1}{3}\kver^{-\flat(j)}\cdot\#\hat{T}\right\} 
\end{align*}
which shows the claim. Furthermore, the crossing happens in one of
the $\hat{S}_{i}$.
\end{proof}
Let us return to our $(n,\flat(n))$ strip $S$. On the bottom of
it, there is a target $(\kver,\flat(n))$ tree $T$, while on top
of it, there is a union of $(\kver,k)$ trees $T'$ with $\pi_{x}(T')\subset\pi_{x}(T)$.
We also recall the parameters $q,m$ and $k$. Let
\begin{align*}
M:=\max\left\{ \flat(m),\,q-1\right\} \qquad k':=\min\left\{ k,\,\flat(m)\right\}  & \,.
\end{align*}
and $\underline{T}$ be a $(\kver,\flat(n))$ tree on the bottom of
$S_{m}$ with $\pi_{x}(\underline{T})=\pi_{x}(T)$. This tree will
act as the target for the survivors of Event $A$. Next, we have to
count the survivors. 

Define $\tilde{T}$ to be the union of $(\kver,q-1)$ trees in $\underline{T}$
satisfying the following: Let $\tilde{T}_{i}$ be a $(\kver,q-1)$
tree inside a $(m,M)$ strip. Then $\tilde{T}_{i}\subset\tilde{T}$
if there are $v_{i}'\in T'$ and $\tilde{v}_{i}\in\tilde{T}_{i}$
such that $v_{i}'\leadsto\tilde{v}_{i}$ inside $S_{m}$. Define the
event
\begin{equation}
\mathfrak{X}:=\Big\{\#\tilde{T}\geq\max\big\{\frac{\kver^{q-2}\cdot\#T'}{8\cdot\kver^{M}},\,\kver^{q-2}\big\}\Big\}\,.\label{eq:Enough}
\end{equation}

\begin{rem*}[On $M,k'$]
 We have to consider $(m,M)$ strips rather than $(m,\flat(m))$
strips because multiple $(m,\flat(m))$ strips might connect to the
same $q-1$ box in the case of $q-1>\flat(m)$. This would result
in double counting for $\tilde{T}$. On the other hand, introducing
$k'$ basically just means that we break up $(\kver,k)$ trees into
smaller $(\kver,k')=(\kver,\flat(m))$ trees so that they act as proper
inputs for the $(m,\flat(m))$ strips.
\end{rem*}
We only count hits of $(\kver,q-1)$ trees since each (single) connection
will yield a full tree after passing through a $q-1$ box (or rather
a $q-1$ column in Event B later).
\begin{lem}[{Probability of ``sufficiently many'' crossings, \cite[Lemma 4.5]{MR2116736}}]
\label{lem:sufficiently-many-crossings}\ Suppose Proposition \ref{prop:Drilling}
holds for $j\leq n-2$. Then
\[
\P(\mathfrak{X})\geq\min\left\{ 0.9,\,\tfrac{1}{8}\kver^{-\flat(m)}\cdot\#T'\right\} \,.
\]
\end{lem}

\begin{proof}
Since $\tilde{T}$ consists of $(\kver,q-1)$ trees and each such
tree has $\kver^{q-2}$ many vertices, we have $\#\tilde{T}\geq\kver^{q-2}$
if and only if $\tilde{T}\neq\emptyset$. In order to show $\#\tilde{T}\geq\kver^{q-2}$,
it therefore suffices to show $T'\leadsto\underline{T}$. The proof
is broken up into cases based on the size of $\#T'$ and the value
of $M$. 
\begin{enumerate}
\item $\#T'\leq8\cdot\kver^{\flat(m)}$ and $M=\flat(m)$. In particular,
$\flat(m)\geq q-1$. Therefore, by Lemma \ref{lem:pooling-strips}
with $S'=S_{m}$, $\hat{S}$ to be a union of $(m,\flat(m))$ strips,
$T^{*}=\underline{T}$ and $\hat{T}=T'$ 
\[
\P(\mathfrak{X})=\P(\#\tilde{T}\geq\kver^{q-2})\geq\P(\exists\text{ a crossing }T'\leadsto\underline{T}\text{ inside }S_{m})\geq\min\Big\{0.9,\,\tfrac{1}{3}\kver^{-\flat(m)}\cdot\#T'\Big\}\,.
\]
\item $\#T'\leq8\cdot\kver^{M}$ and $M=q-1$. Again
\[
\P(\mathfrak{X})=\P(\#\tilde{T}\geq\kver^{M-1})=\P(\#\tilde{T}\geq\kver^{q-2})\,.
\]
Write $T'=\cup_{i=1}^{N}T_{i}$ where each $T_{i}$ is a union of
$(\kappa,k')$ trees in a $(m,\flat(m))$ strip. Then, for all $i$
by Lemma \ref{lem:pooling-strips}
\[
\P\big\{ T_{i}\leadsto\underline{T}\text{ inside a }(m,\flat(m))\text{ strip}\big\}\geq\min\left\{ 0.9,\,\tfrac{1}{3}\kver^{-\flat(m)}\cdot\#T_{i}\right\} \,.
\]
We are done if the minimum for one of the $i$ is $0.9$. Otherwise,
Lemma \ref{lem:pooling-strips} concludes
\[
\P\big\{ T'\leadsto\underline{T}\text{ inside some }(m,\flat(m))\text{ strip}\big\}\geq\min\left\{ 0.9,\,\tfrac{1}{3}\kver^{-\flat(m)}\cdot\#T'\right\} \,.
\]
\item $\#T'>8\cdot\kver^{M}$. This is the case where we actually have to
establish multiple crossings in disjoint regions. Write $T'=\cup_{i=1}^{N'}T_{i}'$
where each $T_{i}'$ is now a union of $k'$ trees that belong to
a union of $(m,M)$ strips $\tilde{S}_{i}$. Do this in a way such
that for each $i$
\[
3\cdot\kver^{M}\leq\#T'_{i}\leq4\cdot\kver^{M}
\]
and such that if $i\neq j$, then the corresponding unions of $(m,M)$
strips $\tilde{S}_{i}$ and $\tilde{S}_{j}$ are disjoint. This is
possible since each $k'$ tree has $\kver^{k'-1}$ vertices and $M\geq\flat(m)\geq k'$.
Thus, $N'$ satisfies
\[
N'\geq\frac{\#T'}{4\cdot\kver^{M}}\geq\frac{8\cdot\kver^{M}}{4\cdot\kver^{M}}=2\,.
\]
By Lemma \ref{lem:pooling-strips}, we have with $\#T_{i}'\geq3\cdot\kver^{M}$
\[
\P\big\{ T_{i}'\leadsto\bar{T}\text{ inside some }(m,M)\text{ strip}\big\}\geq\min\left\{ 0.9,\,\tfrac{1}{3}\kver^{-\flat(m)}\cdot\#T_{i}'\right\} =0.9\,.
\]
Therefore, we have $N'$ independent events with probability greater
or equal to $0.9$. The probability of at least $\lceil N'/2\rceil$
of these happening is $\geq0.9$. Each such event gives us a contribution
of $\kver^{q-2}$ to $\#\tilde{T}$, so we see that under the event
of at least $\lceil N'/2\rceil$ crossings happening
\[
\#\tilde{T}\geq\frac{N'}{2}\cdot\kver^{q-2}\geq\frac{\#T'\cdot\kver^{q-2}}{8\cdot\kver^{M}}\,.
\]
Therefore
\[
\P(\mathfrak{X})\geq\P\Big(\#\tilde{T}\geq\frac{\#T'\cdot\kver^{q-2}}{8\cdot\kver^{M}}\Big)\geq0.9=\min\left\{ 0.9,\,\tfrac{1}{8}\kver^{-\flat(m)}\cdot\#T'\right\} \,.
\]
\end{enumerate}
With this, all cases have been covered.
\end{proof}
This covers event $A$. Next up is event $B$. Take a column of $q-1$
boxes including the $(q,q-1)$ strips inbetween. Let us fix a survivor
$v\in\underline{T}$ from Event $A$, that is, $v$ satisfies $T'\leadsto v$.
We now formalise what is meant by event $B$: 
\begin{defn}[Good $q-1$ columns ]
 Let a column of up to $12\s_{t}$ many $q-1$ boxes be given including
their $(q,q-1)$ strips inbetween. We call it a $q-1$ \textbf{column}
and we call it \textbf{good for} $v,w\in G$ if $v\leadsto w$ inside
$G$.
\end{defn}

\begin{lem}[{Probability of good $q-1$ columns \cite[Lemma 4.6]{MR2116736}}]
\label{lem:Lemma-46}\label{lem:Good-Columns}\ Suppose Lemma \ref{lem:Main-Thm}
holds for $q-1\leq n-2$. Consider a $q-1$ column $G$ and $v,w\in G$
where $v$ is a a vertex on the top and $w$ on the bottom of $G$.
Then,
\[
\P(G\text{ is good for }v,w)\geq0.99\,.
\]
\end{lem}

\begin{proof}
First, we see that $G$ is good for $v$ and $w$ if
\begin{enumerate}
\item all the corresponding $q-1$ boxes and $(q,q-1)$ strips are good
and
\item $v\in\In(\bar{B}_{q-1})$ with $\bar{B}_{q-1}$ being the topmost
$q-1$ box in $G$.
\item $w\in\Out(\underline{B}_{q-1})$ with $\underline{B}_{q-1}$ being
the bottommost $q-1$ box in $G$.
\end{enumerate}
By the induction hypothesis
\[
\P(\text{all of the }q-1\text{ boxes are good})\geq(1-\pp)^{12\s_{t}}\geq1-12\s_{t}\cdot\pp\,,
\]
and
\[
\P(\text{all of the }(q,q-1)\text{ strips are good})\geq(1-\pp)^{12\s_{t}}\geq1-12\s_{t}\cdot\pp\,.
\]
Next, $v\in\In(\bar{B}_{q-1})$ if $v$ lies in good $j$ boxes for
all $j\leq q-1$. The probability of this happening is at least
\[
\P(v\in\In(\bar{B}_{q-1}))\geq1-\sum_{j\geq1}\pp^{j}=\frac{1-2\pp}{1-\pp}\geq1-2\pp\,.
\]
The same holds for $w$. Using $\pp\leq(450\s_{t}\cdot\s_{x})^{-6}$
yields
\[
\P(G\text{ is good for }v,w)\geq1-25\s_{t}\cdot\pp\geq0.99\,,
\]
which finishes the proof. 
\end{proof}
Event $C$ corresponds to Lemma \ref{lem:pooling-strips}. 

\subsection{Drilling: proof of Proposition \ref{prop:Drilling} \label{subsec:Drilling}}

We have gathered all the parts, so it is time to combine them. Unfortunately,
we have to deal with quite a lot of case distinctions.
\begin{proof}[Proof of Proposition \ref{prop:Drilling}]
\ We have already shown the case of $q\leq2$ which also includes
the case of $\min\{m,r\}\leq3$. Now, we may always assume that $m\geq4$
as well as $q\geq3$.We employ our strategy of linking together the
Events $A$, $B$ and $C$, that is, 
\begin{enumerate}
\item[(A)] $\mathfrak{X}$ happens on $S_{1}$. This gives us a collection of
$(\kver,q-1)$ trees $\tilde{T}\subset\bar{T}$ on the bottom of $S_{m}$.
Each such tree has some $v\in\tilde{T}$ with $T'\leadsto v$.
\item[(C)] Consider $T^{*}$ on the top of $S_{r}$ with $\pi_{x}(\tilde{T})=\pi_{x}(T^{*})$.
There exists a crossing of $S_{2}$ intersecting $T^{*}$ and $T$,
i.e., some $T^{*}\ni w\leadsto T$.
\item[(B)] The $q-1$ column of $w\in T^{*},v[w]\in\tilde{T}$ is good.
\end{enumerate}
If all these events hold, then there exists a crossing of $S$ from
$T'$ to $T'$ via
\[
T'\leadsto^{A}\tilde{T}\ni v[w]\leadsto^{B}w\in T^{*}\leadsto^{C}T\,.
\]
By Lemma \ref{lem:sufficiently-many-crossings}
\[
\P(A)=\P(\mathfrak{X})\geq\min\left\{ 0.9,\,\tfrac{1}{8}\kver^{-\flat(m)}\#T'\right\} \,.
\]
Under $\mathfrak{X}$, we have 
\[
\#\tilde{T}\geq\max\left\{ \kver^{q-2},\,\frac{\kver^{q-2}\cdot\#T'}{8\cdot\kver^{M}}\right\} \,.
\]
\begin{itemize}
\item If now $\#T'\leq8\cdot\kver^{M}$ , then $\#T^{*}=\#\tilde{T}\geq\kver^{q-2}$
and by the Lemmas \ref{lem:sufficiently-many-crossings}, \ref{lem:Lemma-46}
\begin{align*}
\P(B,C\,\vert\,A) & \geq\P(\exists\text{ a crossing of }S_{2}\text{ intersecting }T^{*}\text{ and }T\,\vert\,\#T^{*}=\kver^{q-2})\cdot0.99\\
\geq & 0.99\cdot\min\left\{ 0.9,\,\tfrac{1}{3}\kver^{-\flat(r)}\kver^{q-2}\right\} \,.
\end{align*}
If $M=\flat(m)$, then using Equation (\ref{eq:m-r-q-half-n}) from
Lemma \ref{lem:extra-estimates-for-final-lemma} below yields
\begin{align*}
\P(\exists\text{ a cluster in }S & \text{ connecting }T\text{ and }T')\geq\P(A)\cdot\P(B,C\,\vert\,A)\\
\geq & 0.9\cdot\tfrac{1}{8}\kver^{-\flat(m)}\#T'\,\cdot\,0.99\cdot\min\left\{ 0.9,\,\tfrac{1}{3}\kver^{-\flat(r)}\kver^{q-2}\right\} \\
\geq & \#T'\cdot\tfrac{1}{27}\kver^{-\flat(n)+\lfloor q/2\rfloor}\geq\kver^{-\flat(n)}\cdot\#T'\,.
\end{align*}
For the case of $M=q-1$, i.e.\ $8\cdot\kver^{\flat(m)}\leq\#T'\leq8\cdot\kver^{M}$,
using Equation (\ref{eq:m-r-q-half-n}) of Lemma \ref{lem:extra-estimates-for-final-lemma}
and $\flat(m)>m/2\geq\lceil q/2\rceil$ yields
\begin{align*}
\P(\exists & \text{ a cluster in }S\text{ connecting }T\text{ and }T')\geq\P(A)\cdot\P(B,C\,\vert\,A)\\
\geq & 0.9\,\cdot\,0.99\cdot\min\left\{ 0.9,\,\tfrac{1}{3}\kver^{q-2-\flat(r)}\right\} \geq\min\left\{ 0.5,\,\tfrac{1}{4}\frac{\kver^{M-1}\cdot\kver^{\flat(m)}}{\kver^{\flat(m)+\flat(r)}}\right\} \\
\geq & \min\left\{ 0.5,\,\tfrac{1}{4}\frac{\kver^{M-1}\cdot\kver^{\lceil q/2\rceil}}{\kver^{\flat(n)-\lceil q/2\rceil-2}}\right\} \geq\min\left\{ 0.5,\,\frac{\kver^{M}\cdot8}{\kver^{\flat(n)}}\right\} \geq\frac{\#T'}{\kver^{\flat(n)}}\,.
\end{align*}
\item If instead $\#T'\geq8\cdot\kver^{M}$, then using
\[
\#\tilde{T}=\#T^{*}\geq\frac{\#T'\cdot\kver^{q-2}}{8\cdot\kver^{M}}
\]
and Lemma \ref{lem:pooling-strips} and Equation (\ref{eq:m-flatm-q})
gives
\begin{align*}
\P(C\,\vert\,A) & \geq\P\Big\{\exists\text{ a crossing of }S_{2}\text{ intersecting }T^{*}\text{ and }T\,\vert\,\#T^{*}\geq\frac{\#T'\cdot\kver^{q-2}}{8\cdot\kver^{M}}\Big\}\\
 & \geq\min\left\{ 0.9,\,\frac{\#T'\cdot\kver^{q-2}}{8\cdot\kver^{M}}\cdot\frac{1}{3\cdot\kver^{\flat(r)}}\right\} \geq\min\left\{ 0.9,\,\frac{\#T'}{24\cdot\kver^{\flat(n)-1}}\right\} \geq2\frac{\#T'}{\kver^{\flat(n)}}\,,
\end{align*}
where the minimum disappears again from $\#T'\leq\#T\leq\kver^{\flat(n)-1}$.
Lemma \ref{lem:Lemma-46} yields
\[
\P(B\,\vert\,A,C)\geq0.99\,.
\]
Putting everything together, we conclude the $\#T'\geq8\cdot\kver^{M}$
case:
\begin{align*}
\P(\exists\text{ a cluster in } & S\text{ connecting }T\text{ and }T')\geq\P(A)\cdot\P(C\,\vert\,A)\cdot\P(B\,\vert\,A,C)\\
\geq & 0.9\cdot2\cdot\frac{\#T'}{\kver^{\flat(n)}}\cdot0.99\geq\frac{\#T'}{\kver^{\flat(n)}}\,.
\end{align*}
\end{itemize}
This finishes the proof of Proposition \ref{prop:Drilling}.
\end{proof}
\begin{lem}[Extra estimates for final proof]
\label{lem:extra-estimates-for-final-lemma} Let $m\geq4,q\geq3$
and $M=\max(\flat(m),q-1)$. We have
\begin{equation}
\flat(m)+\flat(r)-\lceil q/2\rceil\leq\flat(n)-2\,.\label{eq:m-r-q-half-n}
\end{equation}
Furthermore, we have
\begin{align}
M+\flat(r)-q & \leq\flat(n)-3\,.\label{eq:m-flatm-q}
\end{align}
\end{lem}

\begin{proof}
If $3\leq q\leq8$, then $m+r=n$ by Equation (\ref{eq:m-plus-r-equals-n})
in Lemma \ref{lem:q-estimates}. In particular,
\[
\flat(m)+\flat(r)\leq\flat(n)\implies\flat(m)+\flat(r)-\lceil q/2\rceil\leq\flat(n)-2\,.
\]
If $q\geq9$, then we use
\[
\lceil(2-\cD)^{-1}(\lfloor\cD q\rfloor+1)\rceil\leq q/21+2\leq\lceil q/2\rceil-2
\]
to also obtain Equation (\ref{eq:m-r-q-half-n}) via Equation (\ref{eq:m-plus-r-minus-q-is-n})
in Lemma \ref{lem:q-estimates}
\begin{align*}
m+r-\lfloor\cD q\rfloor & \leq n+1\\
\flat(m)+\flat(r)-\lceil(2-\cD)^{-1}(\lfloor\cD q\rfloor+1)\rceil & \leq\flat(n)\\
\flat(m)+\flat(r)-\lceil q/2\rceil & \leq\flat(n)-2\,.
\end{align*}
For Equation (\ref{eq:m-flatm-q}), we need another case distinction:
If $M=\flat(m)$, then
\begin{align*}
M+\flat(r)-q & =\big\{\flat(m)+\flat(r)-\lfloor q/2\rfloor\big\}-\lceil q/2\rceil\leq\flat(n)-2-1
\end{align*}
Else, we have $M=q-1$, which yields
\begin{align*}
M+\flat(r)-q & =\flat(r)-1=\flat(n)-2-\big\{\flat(m)-\lfloor q/2\rfloor\big\}\,.
\end{align*}
Since $\flat(m)>m/2>\lfloor(m-1)/2\rfloor\geq\lfloor q/2\rfloor$
and $M+\flat(r)-q$ is an integer, this case also implies $M+\flat(r)-q\leq\flat(n)-3$,
i.e., Equation (\ref{eq:m-flatm-q}).
\end{proof}

    \begin{it}
        Acknowledgement.
    \end{it}This work was supported by the German Research Foundation under Germany's Excellence Strategy MATH+: The Berlin Mathematics Research Center, EXC-2046/1 project ID: 390685689, and the Leibniz Association within the Leibniz Junior Research Group on Probabilistic Methods for Dynamic Communication Networks as part of the Leibniz Competition.

\printbibliography[heading=bibintoc]

\end{document}